\documentclass[compress,a4paper,fleqn]{cas-dc}
\usepackage[numbers,sort&compress]{natbib}
\usepackage{graphicx}
\usepackage{subcaption}
\usepackage{tabularx}
\newcolumntype{Y}{>{\centering\arraybackslash}X}
\usepackage{threeparttable}
\usepackage{hyperref}
\hypersetup{hypertex=true,
	colorlinks=true}
\usepackage{cleveref}
\crefname{figure}{Figure}{Figure}
\crefname{algorithm}{Algorithm}{Algorithm}
\usepackage{amsmath,amsfonts}
\usepackage{amssymb}
\usepackage[linesnumbered,lined,ruled,commentsnumbered]{algorithm2e}
\usepackage{array}
\usepackage{caption}
\usepackage{multirow}
\usepackage{array}
\usepackage{booktabs}
\usepackage{ragged2e}

\usepackage{xcolor}
\usepackage{colortbl}
\usepackage{soul}
\soulregister\cite{7}
\soulregister\Cref{7}
\soulregister\ref{7}
\soulregister\eqref{7}

\newcommand{\Low}{\cellcolor{ConstraintLow}L}
\newcommand{\Medium}{\cellcolor{ConstraintMedium}M}
\newcommand{\High}{\cellcolor{ConstraintHigh}\textbf{H}}

\usepackage{tabularx}
\usepackage{array}
\usepackage{booktabs}
\usepackage{ragged2e}
\usepackage{tikz}
\newcommand{\ReviewSolid}{%
	\tikz[baseline=-0.5ex]{%
		\fill (0,0) circle[radius=0.40mm];}}
\newcommand{\ReviewOpen}{%
	\tikz[baseline=-0.5ex]{%
		\draw[line width=0.25pt]
		(0,0) circle[radius=0.40mm-0.125pt];}}

\newcommand{\ReviewThree}{%
	\mbox{\ReviewSolid\,\ReviewSolid\,\ReviewSolid}}
\newcommand{\ReviewTwo}{%
	\mbox{\ReviewSolid\,\ReviewSolid\,\ReviewOpen}}
\newcommand{\ReviewOne}{%
	\mbox{\ReviewSolid\,\ReviewOpen\,\ReviewOpen}}
\newcommand{\ReviewZero}{%
	\mbox{\ReviewOpen\,\ReviewOpen\,\ReviewOpen}}

\definecolor{ConstraintLow}{HTML}{DDEBD5}
\definecolor{ConstraintMedium}{HTML}{FCE5CD}
\definecolor{ConstraintHigh}{HTML}{EAD1D1}

\usepackage{microtype}
\ifXeTeX
  \usepackage{newunicodechar}
  \newunicodechar{–}{--}
  \newunicodechar{’}{'}
\fi

\begin{document}
\let\WriteBookmarks\relax
\renewcommand{\topfraction}{0.90}
\renewcommand{\dbltopfraction}{0.90}
\renewcommand{\textfraction}{0.05}
\renewcommand{\floatpagefraction}{0.80}
\renewcommand{\dblfloatpagefraction}{0.80}
\setcounter{topnumber}{3}
\setcounter{dbltopnumber}{3}
\setcounter{totalnumber}{5}

\title[mode = title]{The Role of Insurance in Power-System Risk Management: A Critical Review}

\author[scu]{Jing Huang}
\author[scu]{Jiale Guo}
\author[scu]{Chenjia Gu}[orcid=0000-0002-9551-7250]
\cormark[1]\ead{gcj0629@scu.edu.cn}
\author[xjtu]{Chengcheng Shao}
\author[scu]{Yue Xiang}[orcid=0000-0001-8456-1195]
\cormark[1]\ead{xiang@scu.edu.cn}
\author[seu]{Tao Qian}
\author[stf]{Alexis Pengfei Zhao}
\author[xjtu]{Xiuli Wang}
\author[scu]{Junyong Liu}
\author[xjtu]{Xifan Wang}
\address[scu]{College of Electrical Engineering, Sichuan University, Chengdu, Sichuan Province 610065, China}
\address[xjtu]{School of Electrical Engineering, Xi'an Jiaotong University, Xi'an, Shaanxi Province 710049, China}
\address[seu]{College of Electrical Engineering, and the Jiangsu Key Laboratory of Smart Grid Technology and Equipment, Southeast University, Nanjing, Jiangsu  Province 210000, China}
\address[stf]{the Department of Energy Science and Engineering, Stanford Doerr School of Sustainability, Stanford University, Stanford, CA 94305, United States}
\cortext[cor1]{Corresponding author.}

\begin{abstract}
	Modern power systems are becoming more renewable-rich, digitalized, and interconnected, making risk propagation and economic consequences increasingly complex. Physical mitigation remains the primary means of risk control, but technical and economic limits leave residual risk. Insurance can complement these measures by transferring associated economic losses and influencing power-system decisions through risk-responsive terms. Despite growing interest, the literature remains fragmented across power-system engineering, actuarial science, and risk finance. This review develops an integrated physical–financial perspective for reviewing insurance in power systems. It examines how power-system risks generate residual economic losses and the conditions under which these losses can be insured and transferred. It then develops an actuarial framework for insured-loss assessment and pricing. Insurance-aware decision-making is synthesized through four coupling levels describing how financial protection interacts with power-system decisions. Further progress is constrained by weak physical–financial loss validation, incomplete representation of systemic and evolving portfolio risk, and limited empirical and actuarial grounding of insurance–decision feedback. Research priorities include stronger physical–financial data linkage, forward-looking portfolio assessment, multi-timescale decision feedback, and empirical validation. Addressing these priorities would help move power-system insurance from fragmented applications toward a more coherent approach to residual-risk management in future power systems.	
\end{abstract}


\maketitle

\section{Introduction}

\subsection{Motivation}

Power systems are rapidly becoming renewable-rich, digitalized, and highly interconnected. This transition is reshaping both the sources of system risk and the mechanisms through which disturbances propagate. Decarbonization increases the penetration of weather-dependent resources and the associated exposure of both supply and demand to common meteorological conditions \cite{staffell2018weather,ref42,gu2026climateAdaptiveReview}, while digitalization strengthens the coupling between information integrity and physical operation \cite{ref121,ref123}. Growing infrastructure interdependencies further create pathways for disturbances to propagate beyond conventional system boundaries \cite{ref65,ref110}. These changes are compounded by extreme weather and climate change, which can alter the frequency, spatial extent, and duration of system stress and introduce increasingly correlated and non-stationary risks \cite{ref21,ref23,ref27}. As a result, disturbances that were once relatively localized or independently treated are increasingly capable of developing into wide-area, long-duration, and compound system disruptions. Power-system risk is therefore becoming more complex, more interconnected, and more difficult to anticipate and manage.

Physical risk mitigation remains the foundation of power-system risk management. A broad range of planning and operational measures, such as equipment hardening, reserve allocation, and emergency restoration, have been developed to reduce the likelihood and consequences of adverse events \cite{ref4,ref7}. These measures are supported by a substantial body of reliability and resilience research on preventing failures, limiting system consequences, and facilitating recovery \cite{ref18,ref22}. Nevertheless, the changing risk environment exposes three limitations that prevent physical measures from providing complete risk coverage.

\textit{i) Challenge of an expanding and evolving risk landscape.}
The scope and characteristics of risks faced by power systems are changing substantially. Beyond conventional equipment contingencies, system operation is increasingly exposed to weather-dependent uncertainty, correlated external disturbances, and interactions with other infrastructures \cite{ref37,ref40}. These risks are also becoming more interdependent, as multiple drivers can occur jointly or sequentially and interact with the prevailing system condition, producing consequences that cannot be adequately characterized by considering individual disturbances in isolation \cite{zscheischler2018compound,zscheischler2020typology,clarfeld2020correlated,sundar2023meteorological}. Climate change and ongoing system transformation further introduce non-stationarity, with both disruptive conditions and system vulnerability evolving over time \cite{turner2019compound,ref29}. Such changes also make risk assessment increasingly dependent on the integration of climate information with power-system modelling across different scales and uncertainty representations \cite{ruan2024interdisciplinary}. Consequently, physical risk-mitigation strategies designed around predefined scenarios and historical regularities may provide incomplete protection against an expanding and evolving risk landscape.

\textit{ii) Challenge of incomplete physical risk coverage.}
Physical risk mitigation is inherently constrained by investment budgets, the response capability of flexible resources, and their availability under disruptive conditions. Expanding protection against increasingly severe or widespread disturbances generally requires greater redundancy and additional flexibility, but such reinforcement is ultimately limited by economic considerations. Even when these resources have been deployed, their ability to sustain system operation remains finite, particularly under prolonged or compound disturbances. Existing studies have shown that coordinated use of local supply, network flexibility, and distributed or mobile resources can substantially reduce service interruptions following major events \cite{ref60,ref58,ref85,ref86,ref88,ref67}. Their effectiveness, however, depends on whether sufficient flexibility remains available and whether the underlying network and supporting resources survive the disturbance. As failures become more widespread or correlated, these conditions become increasingly difficult to maintain. Physical measures can therefore reduce the likelihood and severity of losses, but cannot economically or technically eliminate residual tail risks that may exceed the financial tolerance of utilities, generators, or customers.

\textit{iii) Challenge of residual-loss allocation and risk-reduction incentives.}
Power-system risk is jointly shaped by multiple stakeholders, while the costs of risk reduction and the benefits of avoided losses are often distributed unevenly among them. An actor undertaking risk-reduction investment may therefore bear most of the associated cost but capture only part of the resulting system-wide benefit, weakening the incentive for otherwise valuable actions \cite{ref155}. At the same time, residual losses must be allocated among stakeholders with different preferences, obligations, and capacities to absorb financial losses. Differences in reliability preferences can substantially affect how interruption losses and the value of risk reduction are perceived across stakeholders \cite{ref154}. Related studies further show that the allocation and transfer of outage-related losses can influence both system-level and customer-level decisions \cite{ref203,ref205}. Without a mechanism that assigns a financial value to residual risk and transmits this value to the stakeholders whose decisions influence it, the economic consequences of risk remain only weakly connected to risk-reduction decisions.

Insurance provides a financial complement to the limitations of physical risk mitigation. As the risk environment expands and evolves, insurance requires diverse physical consequences to be translated into stakeholder-specific monetary losses that can be defined within a contractual framework. When physical measures cannot eliminate severe residual losses, insurance and reinsurance provide additional risk-bearing capacity and spread losses beyond a single exposed party \cite{kousky2019disasterInsuranceReview,surminski2016climateResilience}. Insurance can also improve risk-reduction incentives by linking premiums and other contract terms to the risk characteristics and mitigation actions of the insured. Evidence from disaster-insurance studies suggests that such incentives are more effective when mitigation can be observed and reflected in insurance terms \cite{ref174,hudson2016insuranceAdaptation,mol2020moralHazard}. In this way, insurance complements physical risk reduction by transferring residual losses and creating a financial link between risk and risk-reduction decisions.

\begin{figure*}[width=\textwidth,pos=t]
	\centering
	\includegraphics[width=\textwidth]{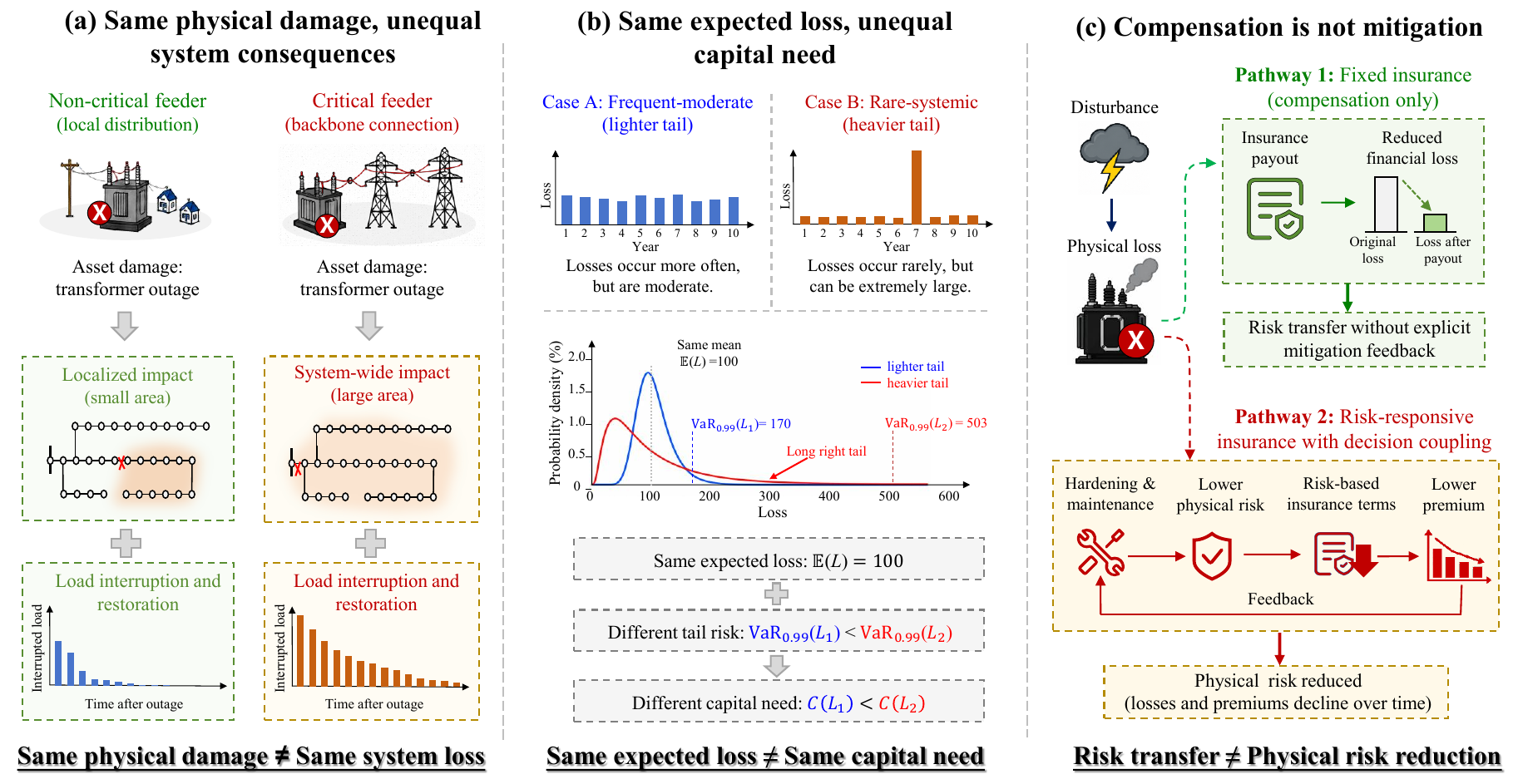}
	\caption{Motivation for integrating insurance into power-system risk management.}
	\label{fig:motivation}
\end{figure*}

Fig.~\ref{fig:motivation} further illustrates why these functions require a power-system-specific perspective. Similar physical asset damage can produce very different system consequences when the affected assets occupy different positions or functions within the network, while similar expected losses can imply markedly different risk-bearing and capital requirements because of differences in tail exposure. Moreover, insurance compensation does not itself reduce physical risk. Its broader risk-governance value emerges when risk-responsive contract terms create feedback to maintenance, hardening, and other mitigation decisions.

Direct research on power-system insurance has addressed an increasingly broad range of risks and applications. Early studies concentrated on comparatively well-bounded operating exposures. Yang et al. linked wind forecast imbalance to insurance settlement \cite{ref150}, while D'Achiardi et al. used reliability contracts to allocate delivery responsibility between renewable and gas producers \cite{ref151}. Subsequent work extended insurance to storage participation \cite{ref152,ref204}, customer reliability differentiation \cite{ref154}, and cyber-physical security \cite{ref127,ref128}. Related cyber-insurance studies also emphasize the governance and information challenges that arise when premiums are expected to reward security effort \cite{ref146,ref147}. More recent studies address catastrophe loss and resilience planning \cite{ref202,ref203}, peer-to-peer outage allocation \cite{ref205}, integrated power--gas risk sharing \cite{xie2023coalitional}, and insurance-integrated operation under extreme weather \cite{wang2026hierarchicalInsurance}. The function of insurance has consequently expanded, yet existing studies remain fragmented across different risks, mechanisms, and applications.

Existing reviews have established strong foundations for different parts of this problem, but they have largely developed along distinct analytical lines. Power-system resilience reviews mainly examine the ability of the system to sustain service under disturbances and restore normal operation afterward \cite{wang2022resilienceReview,younesi2022resilienceReview,shi2022distributionResilienceReview}. More recent studies extend this discussion to microgrid applications and the broader development of resilience research \cite{zidane2025microgridResilience,venkatasubramanian2023resilienceBibliometric}. Climate and extreme-weather reviews instead focus on how changing environmental conditions affect power-system exposure and vulnerability, as well as the need for adaptation \cite{panteli2015weatherResilience,cronin2018climateEnergyReview}. Recent work further treats climate change as a long-term driver of system risk rather than merely a source of individual extreme events \cite{ruan2026longHorizon,bromberger2026extremeWeatherReview}. Reviews of renewable uncertainty and risk-aware planning and operation address how uncertain generation and other sources of variability influence system performance and decision-making \cite{ghanbarzadeh2025capacityPlanningReview,haugen2023uncertaintyReview,wang2025uncertaintyReview}. In parallel, disaster-insurance and catastrophe-risk finance reviews provide a mature body of knowledge on loss assessment, insurance design, and financial risk transfer \cite{kousky2019disasterInsuranceReview,ingels2024climateInsuranceReview,zhao2021catastropheReview}. Together, these research streams cover many of the physical and financial elements relevant to power-system insurance, but they remain largely separated in perspective and methodology.

To address this gap, this review develops an integrated physical–financial risk-management perspective. It focuses on how power-system risks can be translated into insurable losses, actuarially assessed, transferred through insurance mechanisms, and incorporated into power-system decision-making. Table~\ref{tab:positioning} positions this review against the major streams of existing reviews.

\begin{table*}[width=\textwidth,cols=7,pos=t]
	\centering
	
	\caption{Positioning of this review relative to major streams of existing reviews.}
	\label{tab:positioning}
	
	\fontfamily{ptm}\selectfont
	\scriptsize
	\renewcommand{\arraystretch}{1.15}
	\setlength{\tabcolsep}{3pt}
	
	\renewcommand{\tabularxcolumn}[1]{m{#1}}
	
	\begin{tabularx}{\textwidth}{
			@{}
			>{\hsize=1.15\hsize\linewidth=\hsize
				\RaggedRight\arraybackslash}X
			>{\hsize=1.65\hsize\linewidth=\hsize
				\RaggedRight\arraybackslash}X
			>{\hsize=1.40\hsize\linewidth=\hsize
				\RaggedRight\arraybackslash}X
			>{\hsize=0.70\hsize\linewidth=\hsize
				\centering\arraybackslash}X
			>{\hsize=0.70\hsize\linewidth=\hsize
				\centering\arraybackslash}X
			>{\hsize=0.70\hsize\linewidth=\hsize
				\centering\arraybackslash}X
			>{\hsize=0.70\hsize\linewidth=\hsize
				\centering\arraybackslash}X
			@{}}
		\toprule
		
		\textbf{Review stream} &
		\textbf{Core question addressed} &
		\textbf{Dominant perspective} &
		\textbf{Risk and consequences} &
		\textbf{Risk transfer} &
		\textbf{Actuarial assessment} &
		\textbf{Decision coupling} \\
		\midrule
		
		Power-system resilience reviews
		\cite{wang2022resilienceReview,younesi2022resilienceReview,shi2022distributionResilienceReview,zidane2025microgridResilience,venkatasubramanian2023resilienceBibliometric}
		&
		How can power systems withstand, respond to, and recover from disruptive events?
		&
		System-performance oriented: focuses on resilience assessment and system response
		&
		\ReviewThree
		&
		\ReviewOne
		&
		\ReviewZero
		&
		\ReviewZero
		\\[1.2mm]
		
		Climate and extreme-weather reviews
		\cite{panteli2015weatherResilience,cronin2018climateEnergyReview,ruan2026longHorizon,bromberger2026extremeWeatherReview}
		&
		How do weather hazards and climate change alter power-system exposure, vulnerability, and the need for adaptation?
		&
		Hazard- and climate-oriented: focuses on external drivers, evolving exposure, and adaptation
		&
		\ReviewThree
		&
		\ReviewOne
		&
		\ReviewZero
		&
		\ReviewZero
		\\[1.2mm]
		
		Renewable uncertainty and risk-aware planning/operation reviews
		\cite{ghanbarzadeh2025capacityPlanningReview,haugen2023uncertaintyReview,wang2025uncertaintyReview}
		&
		How do uncertain generation and other sources of variability influence system performance and decision-making?
		&
		Decision-oriented: focuses on uncertainty representation and risk-aware decision-making
		&
		\ReviewThree
		&
		\ReviewOne
		&
		\ReviewZero
		&
		\ReviewZero
		\\[1.2mm]
		
		Disaster insurance and catastrophe-risk finance reviews
		\cite{kousky2019disasterInsuranceReview,ingels2024climateInsuranceReview,zhao2021catastropheReview}
		&
		How should uncertain losses be assessed, transferred, pooled, and supported by financial capacity?
		&
		Financial-risk oriented: focuses on loss assessment, insurance design, and financial risk transfer
		&
		\ReviewOne
		&
		\ReviewThree
		&
		\ReviewThree
		&
		\ReviewOne
		\\[1.2mm]
		
		\textbf{This review}
		&
		\textbf{How can power-system risk be translated into insurable loss, actuarially assessed, transferred, and coupled with power-system decisions?}
		&
		\textbf{Integrated physical--financial risk-management perspective}
		&
		\ReviewThree
		&
		\ReviewThree
		&
		\ReviewThree
		&
		\ReviewThree
		\\
		
		\bottomrule
	\end{tabularx}
	
	\par\vspace{0.8mm}
	\raggedright
	\fontfamily{ptm}\selectfont
	\scriptsize
	\textit{Note:}
	Symbols indicate analytical scope.
	\ReviewThree\ = systematic coverage;
	\ReviewTwo\ = substantial coverage;
	\ReviewOne\ = limited coverage;
	\ReviewZero\ = little or no explicit coverage.
	\par
	
\end{table*}

\subsection{Review structure and contributions}

This review makes four main contributions. First, it explains the role of insurance in power-system risk management from the perspective of residual risk, clarifying why physical mitigation alone cannot provide complete protection. Second, it develops an integrated physical--financial framework that connects power-system risk with insurability, loss assessment, risk transfer, and insurance pricing. Third, it examines how insurance interacts with power-system decisions through a four-level coupling framework, moving beyond its conventional role as a post-loss compensation mechanism. Finally, it identifies the main limitations of current research and the methodological advances needed for insurance to become more effectively integrated into power-system risk management.

The review is organized around three guiding questions: how power-system risk becomes an insurable loss and is transferred, how insured losses are assessed and priced, and how insurance interacts with power-system decisions. Table~\ref{tab:guiding-questions} summarizes these questions and their corresponding analytical scope.

\begin{table}[width=\columnwidth,cols=3,pos=t]
	\centering
	\caption{Guiding questions and analytical structure of this review.}
	\label{tab:guiding-questions}
	
	\fontfamily{ptm}\selectfont
	\scriptsize
	\setlength{\tabcolsep}{1.0pt}
	\renewcommand{\arraystretch}{1.05}
	
	\begin{tabularx}{\columnwidth}{
			@{}
			>{\RaggedRight\arraybackslash}p{2.30cm}
			>{\RaggedRight\arraybackslash}X
			>{\RaggedRight\arraybackslash}p{1.55cm}
			@{}}
		
		\toprule
		\textbf{Guiding question} &
		\textbf{How this review addresses it} &
		\textbf{Corresponding section} \\
		\midrule
		
		\textbf{Q1. How does physical power-system risk become an insurable economic loss, and how can that loss be transferred?}
		&
		The review characterizes major risk domains, examines residual loss after physical mitigation, evaluates insurability, and distinguishes the structural dimensions of insurance-enabled risk transfer.
		&
		Section~2: risk, insurability, and transfer structure
		\\
		
		\textbf{Q2. How can insured power-system losses be assessed and translated into insurance prices?}
		&
		The review compares experience-based, system-model-based, and hybrid loss assessment, and examines pure-premium pricing, risk loading, credibility adjustment, and endogenous premium determination.
		&
		Section~3: actuarial assessment and pricing
		\\
		
		\textbf{Q3. How does insurance interact with power-system decisions?}
		&
		The review develops four levels of insurance--power-system decision coupling and examines their implications for planning, operation, and market- and user-side risk allocation.
		&
		Section~4: decision coupling and applications
		\\
		
		\bottomrule
	\end{tabularx}
\end{table}

\begin{figure*}[width=\textwidth,pos=t]
	\centering
	\includegraphics[width=\textwidth]{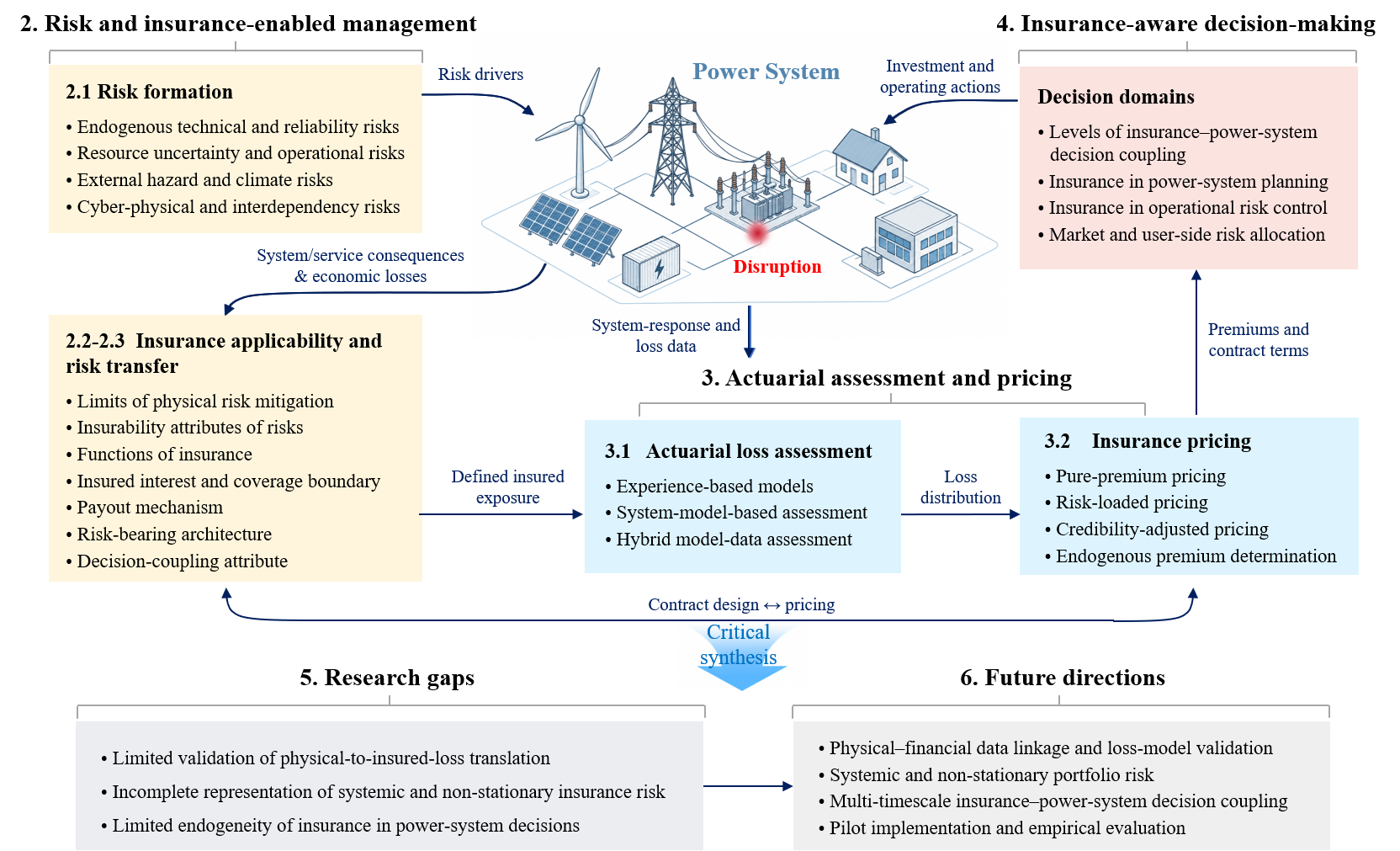}
	\caption{Conceptual framework and review structure for insurance-enabled risk management in power systems.}
	\label{fig:architecture}
\end{figure*}

Fig.~\ref{fig:architecture} presents the overall analytical flow of the review, from physical risk and residual loss to actuarial assessment, risk transfer, and decision feedback. Sections~2--4 develop these three core stages. Section~5 discusses the principal limitations of current research, and Section~6 summarizes the main conclusions and future research directions.


\section{Power-system risk and insurance-enabled management}
\label{sec:risk-transfer}

Insurance enters power-system risk management primarily through the residual risks that remain after physical mitigation. Whether these risks can be effectively transferred depends on how they arise, how their consequences propagate through the power system, and whether the resulting losses satisfy the conditions for insurance. This section first characterizes the major risk sources and the limits of physical mitigation, and then examines the insurability of economic losses arising from residual risk, the role of insurance, and the principal structures through which these risks can be transferred.

\subsection{Risk formation in future power systems}
\label{sec:risk-landscape}

Power-system risk is characterized by the likelihood of adverse system or service outcomes under uncertainty and the severity of the resulting consequences. Its significance depends not only on whether a disturbance occurs, but also on how the power system responds once it occurs. The same initiating event may therefore lead to very different consequences depending on the prevailing system condition and its ability to contain the disturbance and maintain electricity supply \cite{ref2,ref3}. This distinction between an initiating disturbance and its system-level consequence is also reflected in established resilience assessment frameworks \cite{ref18,ref80}.

For analytical purposes, this review groups power-system risks into four domains based on their primary sources and propagation mechanisms. Table~\ref{tab:risk-loss-crosswalk} summarizes how each domain develops from initiating conditions to system- and service-level consequences.

\begin{table*}[width=\textwidth,cols=4,pos=t]
	\centering
	\caption{Formation of major risks in future power systems.}
	\label{tab:risk-loss-crosswalk}
	
	\fontfamily{ptm}\selectfont
	\scriptsize
	\setlength{\tabcolsep}{4pt}
	\renewcommand{\arraystretch}{1.15}
	
	\begin{tabularx}{\textwidth}{
			@{}
			>{\RaggedRight\arraybackslash}p{2.9cm}
			>{\RaggedRight\arraybackslash}p{3.8cm}
			>{\RaggedRight\arraybackslash}p{4.4cm}
			>{\RaggedRight\arraybackslash}X
			@{}}
		\toprule
		
		\textbf{Risk domain} &
		\textbf{Uncertain initiating condition} &
		\textbf{Physical formation mechanism} &
		\textbf{Adverse system or service outcome} \\
		\midrule
		
		Endogenous technical and reliability risks
		&
		Asset deterioration, equipment failure, protection action, or control malfunction
		&
		Failure interaction and cascading propagation shaped by operating state, topology, and corrective control
		&
		Instability, adequacy degradation, or local-to-wide-area interruption of electricity supply
		\\[1.2mm]
		
		Resource uncertainty and operational risks
		&
		Renewable variability, resource unavailability, demand uncertainty, or forecast deviation
		&
		Resource--demand imbalance interacts with available flexibility and network transfer capability
		&
		Curtailment, scarcity, balancing shortfall, or involuntary load loss
		\\[1.2mm]
		
		External hazard and climate risks
		&
		Extreme weather, environmental hazards, and changing climatic conditions
		&
		Hazard, exposure, and vulnerability interact to produce common-cause failures and correlated disruption
		&
		Severe or prolonged supply interruptions, correlated outages, and delayed recovery
		\\[1.2mm]
		
		Cyber-physical and interdependency risks
		&
		Cyber compromise, communication failure, or disruption of dependent infrastructures
		&
		Dependency-mediated propagation removes or constrains observability, control, supply, or restoration capability
		&
		Cross-domain propagation, cascading degradation, or prolonged service interruption
		\\
		
		\bottomrule
	\end{tabularx}
\end{table*}

\subsubsection{Endogenous technical and reliability risks}

Endogenous technical and reliability risks originate from failures within the power system, ranging from equipment outages to failures of protection or control functions. Conventional adequacy and security assessments mainly represent such risks through component availability and predefined contingencies \cite{ref10}. This approach is well suited to routine failures, but becomes less complete when the same outage can produce very different consequences under different system conditions. 
More recent studies therefore place greater emphasis on how the operating state and subsequent system response shape the consequences of component failures \cite{ref81,ref197}. The focus is consequently shifting from isolated contingencies toward correlated and cascading failures.

Once an initial failure alters the operating condition of the network, it may increase the likelihood or severity of subsequent outages. Clarfeld et al. show that correlation among initiating outages can materially change cascading-blackout risk even when the network structure itself is unchanged \cite{clarfeld2020correlated}. Similar effects arise in distribution systems, where the consequences of local failures depend strongly on network configuration and the ability to restore service \cite{ref8,ref15}. Under stressed conditions, common-cause or cascading outages can further reduce the set of feasible operating responses \cite{ref82}. Endogenous technical risk therefore depends increasingly on the interaction between component failures and the system state in which they occur.

\subsubsection{Resource uncertainty and operational risks}

Resource uncertainty arises from weather-driven variability in supply and demand and becomes operationally significant when the resulting imbalance exceeds available balancing capability. Short-term forecast errors affect reserve requirements and dispatch decisions \cite{staffell2018weather,ref184,ref185}, whereas persistent low-output periods create longer-duration scarcity exposure that may exceed the capability of storage, demand response, or network support \cite{ref45,ref46}. Recent studies further show that renewable droughts can be compound and climate conditioned, making resource adequacy increasingly dependent on the duration and persistence of adverse conditions \cite{bracken2025futureDrought,bracken2025seasonalDrought}. The key risk therefore lies in the mismatch between uncertain source--load conditions and the flexibility that can be delivered when and where it is required.

This mismatch is intensified by temporal persistence and spatial
correlation. Prolonged renewable shortfalls increase the duration
for which flexibility must be sustained \cite{ref45,ref46}.
Spatial diversification can smooth aggregate renewable output,
but its benefit weakens when neighbouring regions are exposed
to the same persistent weather regimes
\cite{grams2017spatialWind}. Compound weather conditions can further create simultaneous supply and demand stress. Recent ``power drought'' studies explicitly represent such source--load coupling \cite{bian2025powerDrought}. Sundar et al. similarly show that large-scale weather patterns can jointly increase demand and suppress renewable supply, creating periods of concentrated operational risk \cite{sundar2023meteorological}. Resource uncertainty therefore becomes most consequential when persistent and correlated source--load deviations coincide with insufficient or inaccessible flexibility.

\subsubsection{External hazard and climate risks}

\begin{figure}[!htbp]
    \centering
    \includegraphics[width=\columnwidth]{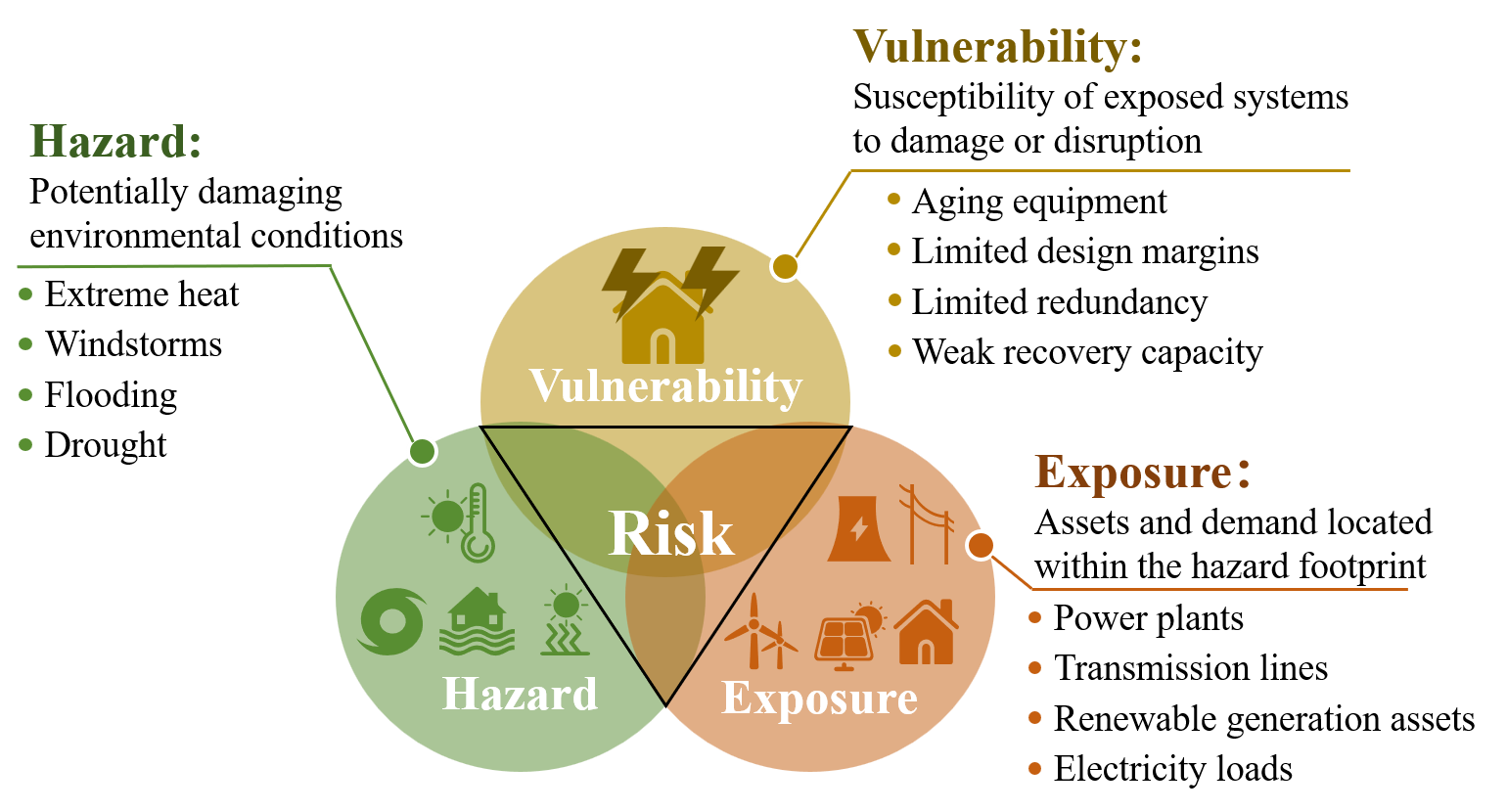}
    \caption{Conceptual representation of risk as the interaction among hazard, exposure, and vulnerability.}
    \label{fig:hev-risk}
\end{figure}

External hazard and climate risks arise from the interaction of hazard, exposure, and vulnerability, as illustrated in Fig.~\ref{fig:hev-risk}. Hazard describes the potentially damaging environmental condition, exposure identifies the assets, resources, and demand subject to that condition, and vulnerability determines how severely the exposed system is affected \cite{chen2021ipccrisk}. In power systems, these three elements are not independent: the same weather event can produce very different consequences depending on where critical assets are located, how system loading and resource availability evolve, and how strongly exposed components respond to the disturbance. Climate change further makes this risk structure increasingly non-stationary by altering both the characteristics of hazardous conditions and the physical system exposed to them \cite{gu2026climateAdaptiveReview}. Changes in temperature and water availability can affect conventional generation capability \cite{ref25,ref30}, while climate-sensitive demand and renewable resources reshape the timing and magnitude of system stress \cite{ref32,ref37}. Rising temperatures can simultaneously reduce transmission
capability and increase peak demand \cite{ref34,ref38}, while asset
vulnerability can evolve as aging and degradation interact with
extreme-weather exposure \cite{ref51}. Historical event frequencies and fragility relationships may therefore become progressively less representative of future risk, motivating climate-resilience studies that explicitly consider evolving hazard baselines and system vulnerability over long planning horizons \cite{xu2024renewableClimate,ruan2026longHorizon}.

A second defining feature is the dependence introduced by common environmental conditions. The same weather system can affect multiple parts of the power system simultaneously, creating correlated consequences that are not well represented by independent contingencies. Turner et al. show that joint climate effects on electricity demand and hydropower can produce shortfalls that are missed when the two pathways are assessed separately \cite{turner2019compound}. Fragility-based weather models likewise show how a common hazard field can generate spatially correlated component failures \cite{ref21,ref52}. St{\"u}rmer et al. further demonstrate that the system consequence of storm damage depends strongly on where failures occur and how they propagate through the network \cite{sturmer2024texas}. This broader dependence across hazards and locations is also emphasized in compound-event research \cite{zscheischler2020typology,raymond2020connected,ref54}. External hazard and climate risk is therefore shaped by event severity and by the interaction between a common hazard and the exposed system.

\subsubsection{Cyber-physical and interdependency risks}

Cyber-physical and interdependency risks arise because power-system operation increasingly depends on digital systems and other infrastructures \cite{ref135}. Their defining feature is dependency-mediated propagation: a disruption originating outside the electrical network can remove functions required for normal operation or recovery, thereby translating an external failure into a power-system consequence. Studies of coupled electricity--gas systems show that such dependencies can materially change feasible operating conditions during disruptions \cite{ref65}, while broader interdependent-infrastructure models demonstrate that cross-system propagation can also alter restoration priorities and outcomes \cite{ref110,ref200}. The resulting risk therefore depends not only on the reliability of individual systems, but also on the dependencies through which disturbances are transmitted to the power system.

Cyber risk represents a particularly direct form of this propagation because disruptions in the information layer can alter physical system behaviour. Manipulated measurements can compromise state estimation \cite{ref122,ref139}, while attacks on control functions may propagate through closed-loop interactions and affect physical operation \cite{ref124,ref125}. As digital technologies become more deeply integrated into power-system operation, these cyber--physical pathways create increasingly important sources of common vulnerability \cite{ref16,ref17}. The central concern is therefore not simply whether a cyber or external-infrastructure disruption occurs, but how far its effects can propagate through the coupled system. Although individual attack paths and specific infrastructure dependencies have been studied extensively, the interaction of multiple dependencies remains less systematically understood.

Across the four domains, power-system risks follow a common progression from an initiating disturbance or uncertainty to system response, service disruption, and ultimately economic consequence. This shared structure provides a consistent basis for comparing otherwise heterogeneous risks. The mechanisms governing this progression, however, differ across domains and lead to distinct loss patterns. These differences need to be reflected in how the resulting risks are characterized and managed.

\subsection{Applicability and functions of insurance}
\label{sec:risk-transfer-architecture}

\subsubsection{Limits of physical risk mitigation}

Physical risk mitigation remains the primary means of reducing the likelihood and consequences of power-system disruptions. It encompasses preventive measures such as reinforcement and topology design \cite{ref61,ref73}, as well as corrective and recovery measures including network reconfiguration and restoration \cite{ref107,ref88}. These measures can substantially reduce physical risk, but the protection they provide is bounded by technical, resource, and economic constraints.

\textit{1) Technical and resource limits.}
Physical mitigation is bounded by system capability and by the resources that remain available during a disruption. Measures designed for a given range of system stress can lose effectiveness once that range is exceeded. Existing studies show that islanding and local supply can preserve service after major disturbances, but only when sufficient resources remain available within the affected area \cite{ref60,ref57,ref58}. Restoration measures face a similar constraint because their effectiveness depends on the post-event availability and accessibility of recovery resources \cite{ref75,ref94}. Physical mitigation can therefore reduce risk within a feasible protection range, but cannot provide complete coverage against all possible system conditions.

\textit{2) Economic limits.}
Even when additional protection is technically feasible, extending physical mitigation toward more extreme or less probable events may not be economically justified. Risk-informed planning therefore requires balancing the reduction in expected disruption against the additional investment needed to achieve it. Existing studies show that the preferred level and composition of protection depend strongly on investment cost, available budget, and the value assigned to risk reduction \cite{ghosh2023resilienceInvestment,poudyal2026resiliencePlanning}. Cost--benefit analysis further shows that additional protection is justified only when its incremental benefit outweighs the associated cost \cite{ciapessoni2023costBenefit}. As physical protection is extended toward increasingly remote loss scenarios, this condition becomes harder to satisfy. Complete elimination of loss is therefore neither a realistic nor an economically efficient objective.

The risk remaining after technically feasible and economically justified mitigation constitutes residual risk. Although physical measures can reduce system exposure and consequences, they cannot eliminate the possibility of economic loss. Insurance provides a means of transferring part of the economic losses arising from this residual risk to external risk-bearing capacity. Disaster-risk research has long recognized financial protection as an important complement to physical prevention and recovery measures \cite{kousky2019disasterInsuranceReview,surminski2016climateResilience}, and emerging power-system studies have begun to examine the joint allocation of resources between physical mitigation and financial risk transfer \cite{ref203}. Insurance therefore provides an additional layer of financial protection against losses that cannot be fully prevented through physical mitigation.

\subsubsection{Insurability attributes of power-system risks}

Classical insurance theory evaluates insurability through conditions concerning the randomness, definition, estimation, diversification, and economic transfer of losses \cite{berliner1982limits,berliner1985large}. Subsequent studies have extended this discussion to emerging risks such as climate and cyber risk, where strong dependence, limited data, and changing risk conditions pose additional challenges to underwriting \cite{charpentier2008insurability,biener2015cyber}. For power-system applications, this review consolidates these considerations into six attributes: fortuity, identifiability, measurability, predictability, diversifiability, and economic insurability. Together, they provide a common basis for assessing whether power-system consequences can be converted into insurable losses and transferred through insurance.

\textit{1) Fortuity.}
Insurance requires the occurrence and magnitude of a covered loss to retain a meaningful degree of randomness and to remain outside the deliberate control of the insured \cite{berliner1982limits,berliner1985large,biener2015cyber}. Many power-system disturbances satisfy this condition at their origin, including forced equipment outages, weather-driven failures, and unexpected resource or demand deviations. The distinction becomes less straightforward because system operators actively modify the consequences of these disturbances through preventive and corrective actions. For example, deliberate de-energization may create an interruption while reducing exposure to a larger wildfire loss \cite{vazquez2022wildfire,ref117}. Fortuity in power systems concerns the uncertainty of the initiating disturbance and the distinction between covered losses and intentional response actions.

\textit{2) Identifiability and verifiability.}
An insurable event must be clearly defined so that its occurrence and resulting loss can be identified and independently verified. Power systems provide extensive operational measurements that can support this requirement. Smart-meter and network data, for example, can help detect and locate line outages \cite{liao2022outage}. Verification becomes more difficult when losses emerge through cascading processes or interactions across coupled systems, because the observed consequence may not be attributable to a single initiating event. Similar challenges arise in cyber risk, where an incident may be observable even when the underlying cause or pathway remains uncertain \cite{ref198,he2024cyberReview}. Insurability therefore requires a clear event boundary that links the verified disturbance to the system and service losses covered by the contract.

\textit{3) Measurability and valuation.}
Once a loss event has been identified, its consequence must be quantified and assigned an economic value. Power-system models can readily measure physical and service outcomes, but translating these outcomes into monetary loss is less straightforward. The economic value of the same service interruption can differ substantially across affected users and operating conditions. Ovaere et al. show that more detailed value of lost load (VoLL) estimates can materially alter reliability decisions \cite{ovaere2019voll}, while Gorman emphasizes that VoLL should be matched to the specific decision context and affected customer group \cite{gorman2022vollReview}. Insurance valuation therefore requires a clear definition of the loss being covered and the party to whom that loss is attributed. The same physical event may generate different economic claims depending on the contractual and market context.

\textit{4) Predictability.}
Insurance pricing requires loss frequency and severity to be estimated with sufficient reliability. Historical data can provide a reasonable basis for routine failures, but are often sparse for rare or emerging risks. Catastrophe modelling therefore uses physically plausible simulations to supplement limited observations \cite{steptoe2022catModeling}, while cyber-insurance research highlights the continued scarcity and inconsistency of usable claims data \cite{eling2022cyberSeverity}. Predictability is further challenged when the underlying risk environment changes over time. Growing weather dependence, shifts in resource adequacy, and evolving exposure to extreme events can make historical loss experience less representative of future conditions \cite{staffell2018weather,sundar2023meteorological,panteli2015weatherResilience}. The key requirement is therefore that future loss distributions can be estimated and updated with sufficient confidence over the period for which insurance coverage is provided.

\textit{5) Diversifiability.}
Risk pooling depends on whether losses across the insured portfolio can be sufficiently diversified. This condition is challenging in power systems because common disturbances and network interactions can create strong dependence among otherwise separate risks. Correlated initiating failures, for example, can materially increase blackout probability and severity \cite{clarfeld2020correlated}, while catastrophe-insurance research similarly shows that dependence can substantially reduce the diversification available to insurers \cite{wu2020catInsuranceEquilibrium}. Cyber insurance provides a related example, where common vulnerabilities and contagion can transform individually manageable risks into systematic or systemic exposure \cite{awiszus2023cyberPricing}. Portfolio-level insurability therefore depends on both the expected loss and the degree of loss dependence within the insured pool.

\textit{6) Economic insurability.}
Meeting the preceding conditions does not necessarily make a risk economically insurable. The premium must cover expected claims and the additional cost of uncertainty and capital, while remaining acceptable to the insured \cite{berliner1982limits,berliner1985large,charpentier2008insurability}. This creates an economic boundary between risks that can be quantified and those that can be transferred at a sustainable price. Existing power-system studies illustrate this transition by translating reliability or outage consequences into monetary losses and incorporating them into premium and risk-pooling calculations \cite{yang2020premium,lau2023mutual}. Evidence from disaster insurance further shows that the cost of transferring risk rises as coverage extends toward more severe loss layers \cite{cebotari2020sovereignInsurance}. Practical coverage may therefore remain limited when tail losses are large, dependence is strong, or uncertainty is high.

Table~\ref{tab:insurability-attributes} summarizes how these six attributes apply across the four power-system risk domains. Their importance varies with the underlying risk and the way losses are formed, observed, and transferred. Insurability can therefore be understood as a property of the resulting loss and its transfer arrangement.

\begin{table*}[width=\textwidth,cols=7,pos=t]
	\centering
	\caption{Classical insurability criteria and their implications for power-system risks.}
	\label{tab:insurability-attributes}
	
	\fontfamily{ptm}\selectfont
	\scriptsize
	\setlength{\tabcolsep}{4pt}
	\renewcommand{\arraystretch}{1.15}
	
	\begin{tabularx}{\textwidth}{
			@{}
			>{\RaggedRight\arraybackslash}p{2.8cm}
			>{\RaggedRight\arraybackslash}p{3.7cm}
			>{\RaggedRight\arraybackslash}X
			>{\centering\arraybackslash}p{0.72cm}
			>{\centering\arraybackslash}p{0.72cm}
			>{\centering\arraybackslash}p{0.72cm}
			>{\centering\arraybackslash}p{0.72cm}
			@{}}
		\toprule
		
		\multirow{2}{=}{\textbf{Insurability criterion}} &
		\multirow{2}{=}{\textbf{Power-system implication}} &
		\multirow{2}{=}{\textbf{Key constraints}} &
		\multicolumn{4}{c}{\textbf{Constraint intensity across risk domains}} \\
		\cmidrule(lr){4-7}
		
		& & &
		\textbf{D1} &
		\textbf{D2} &
		\textbf{D3} &
		\textbf{D4} \\
		\midrule
		
		Fortuity
		&
		Uncertain and non-deliberate loss
		&
		Preventive actions; operator intervention
		&
		\Low & \Medium & \Low & \Medium
		\\
		
		Identifiability and verifiability
		&
		Event detection and consequence reconstruction
		&
		Cascading sequences; topology/control changes; cross-infrastructure propagation
		&
		\Medium & \Medium & \Medium & \High
		\\
		
		Measurability and valuation
		&
		Translation of physical and service consequences into economic loss
		&
		Multiple loss channels; stakeholder heterogeneity; double counting
		&
		\Medium & \High & \High & \High
		\\
		
		Predictability
		&
		Estimation of loss frequency and severity
		&
		Sparse extremes; system evolution; non-stationarity; evolving threats
		&
		\Medium & \High & \High & \High
		\\
		
		Diversifiability
		&
		Pooling without excessive simultaneous loss
		&
		Common-cause failures; spatial correlation; cascading; shared vulnerabilities
		&
		\High & \Medium & \High & \High
		\\
		
		Economic insurability
		&
		Transfer at an economically viable price
		&
		Tail severity; accumulation; model uncertainty; capital requirement
		&
		\Medium & \Medium & \High & \High
		\\
		
		\bottomrule
	\end{tabularx}
	
	\par\vspace{1mm}
	\raggedright
	\fontfamily{ptm}\selectfont
	\scriptsize
	\textit{Note:}
	D1 denotes endogenous technical and reliability risks,
	D2 denotes resource uncertainty and operational risks,
	D3 denotes external hazard and climate risks,
	and D4 denotes cyber-physical and interdependency risks.
	H, M, and L indicate high, moderate, and low levels of constraint, respectively.
	\par
	
\end{table*}

\subsubsection{Functions of insurance in power-system risk management}

\begin{figure}[!htbp]
    \centering
    \includegraphics[width=\columnwidth]{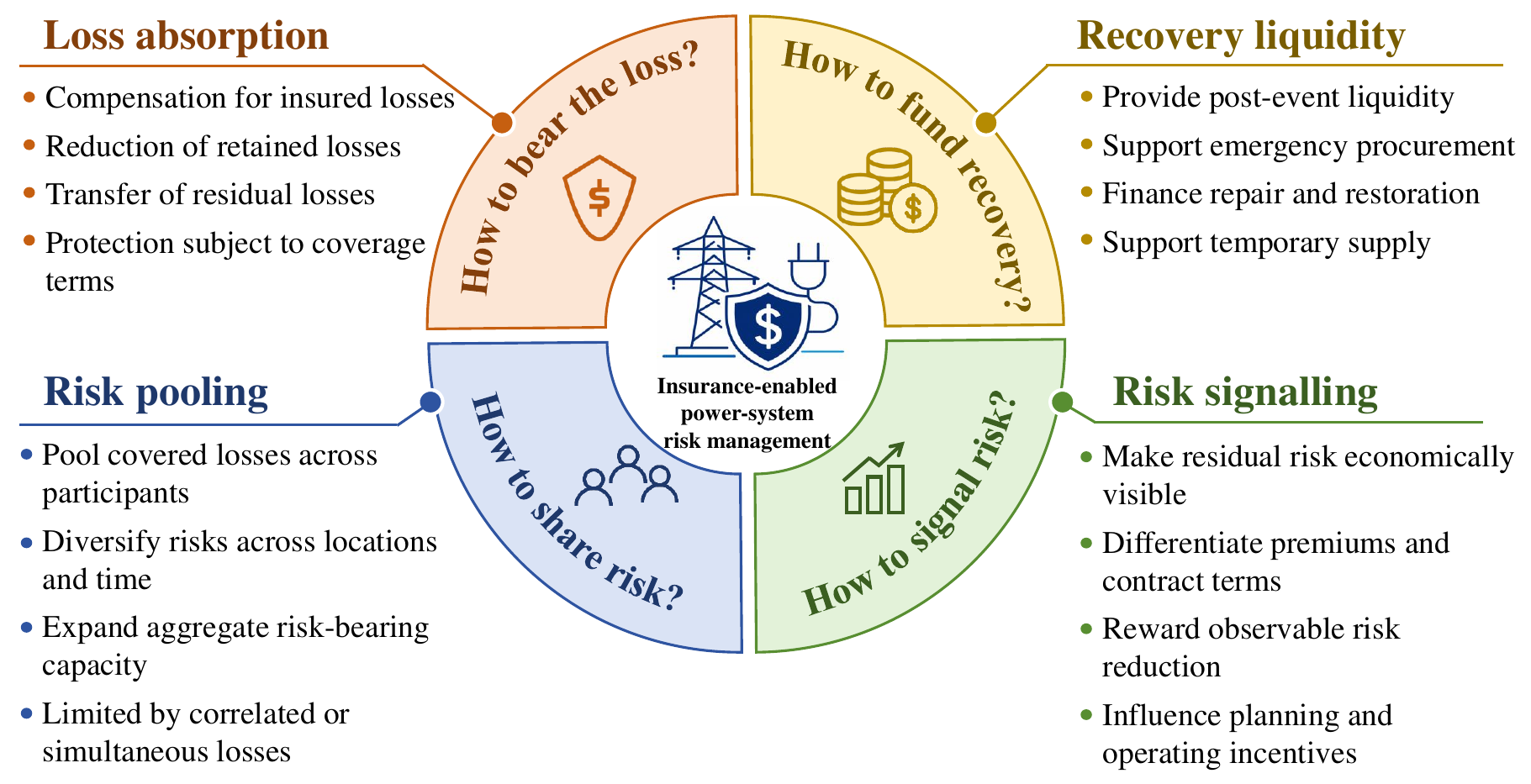}
    \caption{Functions of insurance in power-system risk governance.}
    \label{fig:insurance-functions}
\end{figure}

Once residual power-system losses satisfy the relevant conditions for insurance, insurance can complement physical risk mitigation through several financial functions. As summarized in Fig.~\ref{fig:insurance-functions}, these functions are grouped into four categories: loss absorption, recovery liquidity, risk pooling, and risk signalling.

\textit{1) Loss absorption.}
Loss absorption refers to the use of insurance capital to compensate part of the economic loss that remains after physical mitigation and system response. The covered loss depends on the protected interest and may differ across stakeholders. Unlike physical measures that reduce the underlying system consequence, insurance reduces the portion of economic loss ultimately retained by the affected party. Existing catastrophe and power-system planning studies show that such coverage can substantially reduce retained losses, particularly under severe events \cite{ref202,ref203}. The extent of this protection depends on the design of the insurance contract and the share of loss retained by the insured.

\textit{2) Recovery liquidity.}
Insurance can provide post-event liquidity that supports repair and restoration when financing is itself a constraint. Its value is therefore greatest when recovery resources are physically available but cannot be mobilized quickly because of limited funding. Disaster-risk studies recognize this financing role as an important complement to physical prevention \cite{kousky2019disasterInsuranceReview}. Power-system restoration studies likewise show that recovery depends on the timely mobilization and coordination of available resources \cite{ref109,ref111}. Insurance can accelerate this process when financial capacity is binding, but has limited effect when recovery is constrained primarily by physical availability. Existing catastrophe-insurance models have begun to link insurance settlement with restoration expenditure \cite{ref202}, although the effect of financing on actual recovery trajectories remains only partially represented.

\textit{3) Risk pooling.}
Risk pooling spreads covered losses across a broader risk-bearing base rather than leaving the full consequence with a single exposed party. Existing studies have explored this principle through mutual, peer-to-peer, and cooperative insurance arrangements in energy and infrastructure settings \cite{lau2023mutual,ref205,xie2023coalitional}. The benefit of pooling, however, depends on the extent to which losses can be diversified across participants. Strong dependence among insured risks can produce simultaneous claims and reduce the effectiveness of diversification. 
Therefore, risk pooling can reduce individual exposure and improve collective loss sharing, but its effectiveness depends on the extent to which losses remain diversifiable across participants.

\textit{4) Risk signalling.}
Insurance can also influence power-system decisions by making the cost of residual risk visible to the insured. When risk-reduction actions lead to a credible decrease in expected insured loss, this change can be reflected in contract terms and thereby alter the economic incentive for mitigation. Evidence from natural-hazard insurance shows that risk-based pricing and premium discounts can encourage protective investment when mitigation is observable and recognized by the insurer \cite{hudson2016insuranceAdaptation,mol2020moralHazard}. Similar opportunities exist in power systems because maintenance condition and other operational information can provide observable signals of risk-reduction effort \cite{ref112,ref114,ref115}. Effective signalling therefore requires a credible and verifiable link between the action taken and the resulting change in insured loss.

\subsection{Structural dimensions of power-system insurance}
\label{sec:configuration-logic}

\begin{figure}[!htbp]
    \centering
    \includegraphics[width=\columnwidth]{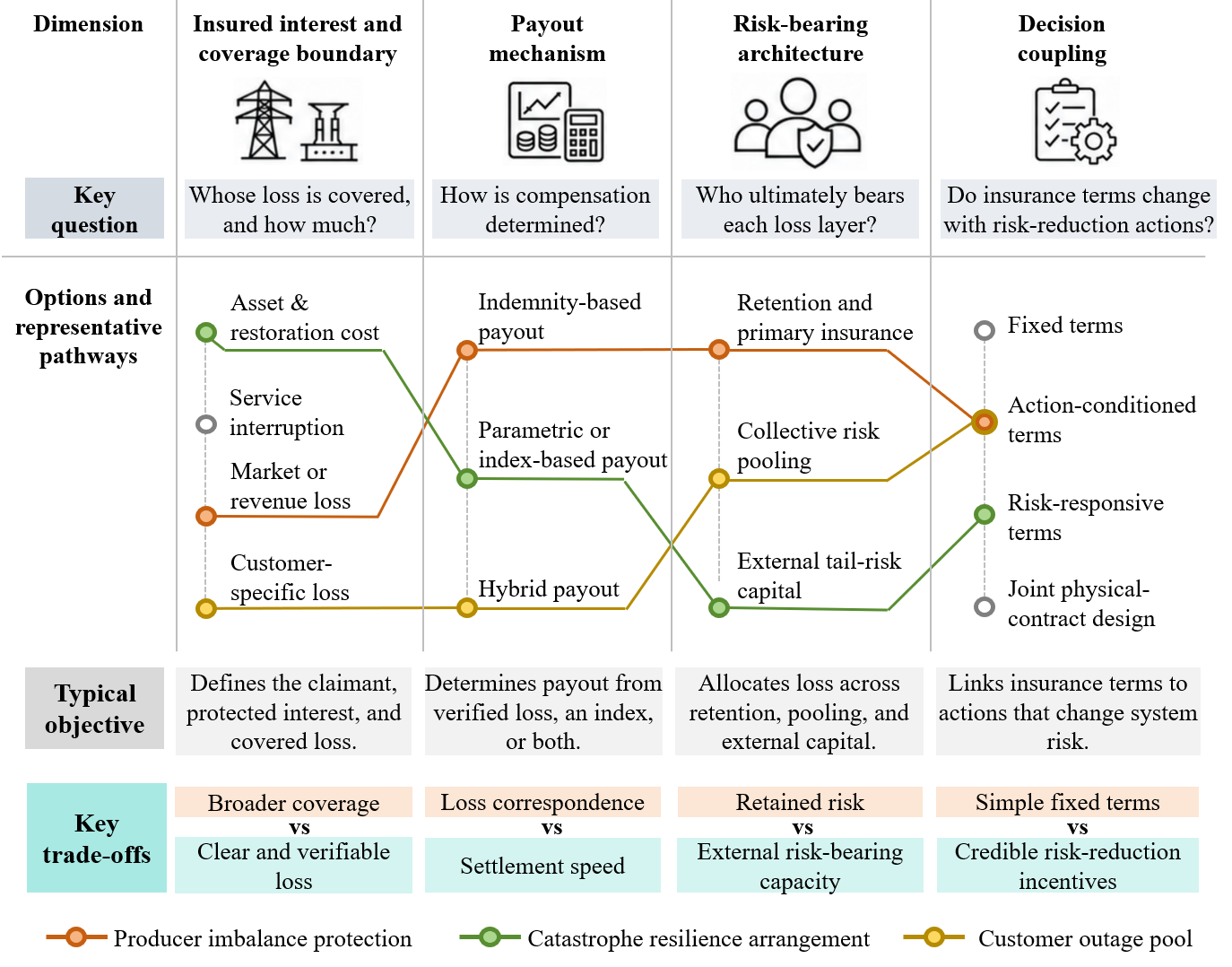}
    \caption{Four structural dimensions and representative configuration pathways of power-system insurance.}
    \label{fig:contractconfiguration}
\end{figure}

The effectiveness of insurance depends on both insurability and contract structure. For power-system applications, this review characterizes insurance arrangements along four dimensions, which provide a common basis for comparing different insurance designs. Fig.~\ref{fig:contractconfiguration} illustrates how representative arrangements combine choices across these dimensions to form different pathways of risk transfer.

\subsubsection{Insured interest and coverage boundary}

The insured interest and coverage boundary determine whose economic loss is protected and which losses fall within the contract. This distinction is important in power systems because the same physical disturbance can create different economic consequences for different stakeholders. Existing studies reflect this diversity by covering producer-side losses in renewable-energy applications \cite{ref150,liao2022solarIndex}, utility losses associated with major disruptions \cite{ref202,ref203}, and customer interruption losses under differentiated reliability preferences \cite{ref154,ref155,ref205}. The broader VoLL literature likewise shows that the economic consequence of an interruption depends strongly on the affected customer and the context in which the loss is valued \cite{gorman2022vollReview}. A common physical event therefore does not correspond to a single economic loss.

Insurance contracts must consequently define the protected interest and the portion of loss transferred to the insurer. This boundary also distinguishes insurance from mechanisms that reward physical availability or delivery performance. Reliability contracts impose obligations during stressed system conditions \cite{ref151}, while capacity mechanisms and reliability options provide compensation according to predefined market rules \cite{ref165,ref167}. Insurance instead transfers a specified economic loss associated with impairment of the insured interest. An important unresolved issue is how to treat economically linked or overlapping claims when the same disturbance creates losses for multiple stakeholders.

\subsubsection{Payout mechanism}

The payout mechanism determines how compensation is calculated once the insured interest has been defined. Three representative structures are considered here: indemnity-based, parametric or index-based, and hybrid. Their main difference lies in how closely the payout follows the realized economic loss.

\textit{1) Indemnity-based payout.}
Indemnity insurance bases compensation on the verified loss incurred by the insured party. It therefore provides a direct link between payout and realized economic consequence, provided that the covered loss can be reliably assessed. Existing power-system studies derive insured losses from simulated asset damage and service interruption \cite{ref202}, or from the economic exposure created by renewable forecast errors and imbalance settlement \cite{ref150}. The main challenge is that power-system losses are often shaped by the system state in which the disturbance occurs, making event attribution and counterfactual loss assessment difficult. Greater correspondence with realized loss may therefore require more detailed modelling and claims verification.

\textit{2) Parametric or index-based payout.}
Parametric insurance determines compensation from a pre-specified observable index rather than from complete ex-post loss assessment. Renewable-energy applications have used weather and production indices to represent generation-related exposure \cite{liao2022solarIndex,mosquera2022weather,han2019windIndex}. Because the trigger and payout rule are defined in advance, settlement can be relatively transparent and rapid. The main limitation is basis risk, which arises when the index-based payout does not match the actual economic loss \cite{clement2018basisRisk,kousky2021parametric}. This mismatch can be particularly important in power systems because the same external condition may lead to very different losses under different system states. More granular data and improved index design can reduce basis risk
\cite{niakh2025basis}. Recent utility applications further show that composite
indices combining hydrological, demand-related, and market signals can better track
the underlying financial exposure than separate single-index products
\cite{amonkar2025compositeIndex}, although some mismatch generally remains.

\textit{3) Hybrid payout.}
Hybrid arrangements combine loss-based and index-based information so that rapid observable triggers can be used without relying entirely on a proxy for final loss. Existing studies have proposed blended structures that switch or combine indemnity and index components according to the quality of the available loss signal \cite{zhang2024blended,lopez2026indexHybrid}. Such designs remain relatively limited in power-system applications, but are potentially useful where an event can be detected quickly while its final economic consequence requires more detailed assessment.

The appropriate payout mechanism therefore depends on the relationship between observable information and the insured loss. Indemnity-based arrangements emphasize correspondence with realized loss,  parametric or index-based arrangements emphasize speed and verifiability, while hybrid structures seek a practical compromise between the two. For power-system insurance, the central issue is whether the selected payout rule provides a sufficiently reliable representation of the economic loss that the contract is intended to transfer.

\subsubsection{Risk-bearing architecture}

Risk-bearing architecture determines where transferred losses ultimately reside and how much financial capacity is available to absorb them. Three representative layers are considered:

\textit{1) Retention and primary insurance.}
At the most basic level, the insured retains part of the loss and transfers the covered portion to a primary insurer. Most existing power-system studies adopt this bilateral structure because it is sufficient for analysing coverage demand and insurance--decision interactions. The insurer's own capacity is usually treated as given, however, so the ability to absorb large or correlated claims is rarely examined explicitly.

\textit{2) Collective risk pooling.}
Collective arrangements spread losses across a broader group of participants. Existing studies have explored mutual, coalitional, and peer-to-peer structures for sharing heterogeneous power-system risks \cite{ref128,xie2023coalitional,ref205}. Their effectiveness depends on the extent to which losses can be diversified across the pool. When participants are exposed to common drivers, aggregate claims may become highly correlated and substantially reduce the benefit of pooling. Catastrophe-insurance research therefore emphasizes the joint role of diversification and available capital in determining risk-bearing capacity \cite{dionne2022catCapacity}.

\textit{3) External tail-risk capital.}
Losses that exceed the capacity of the primary insurer or pool can be transferred further to external sources of capital. Reinsurance and catastrophe bonds are widely used to provide additional capacity for severe loss layers \cite{ref175,lakdawalla2012catBondReinsurance,nematshahi2025wildfireCatBond}, while broader disaster-finance frameworks combine insurance with other financing instruments to support extreme losses \cite{clarke2017sovereignRiskFinance}. Catastrophe-market studies likewise show that the amount of risk that can ultimately be transferred depends on the available capital base and the dependence of aggregate claims \cite{dionne2022catCapacity}. These external layers remain underdeveloped in power-system insurance models.
Recent work has begun to explicitly introduce reinsurance into transmission-system
catastrophe-risk management, linking heavy-tailed loss assessment with insurer
retention and solvency considerations \cite{zhao2026multiStakeholder}. Nevertheless,
most existing studies still focus more on allocating losses than on testing whether
sufficient risk-bearing capacity remains available under severe correlated events.

\subsubsection{Decision coupling}

Decision coupling describes the extent to which insurance terms respond to decisions that change power-system risk. In weakly coupled arrangements, insurance remains largely external to power-system decision-making. Decisions may affect the realized loss, but the contract itself does not adjust in response. Stronger coupling emerges when observable risk-reduction actions are reflected in insurance conditions, creating a financial incentive for mitigation. Existing studies have begun to incorporate such risk-responsive mechanisms by linking cyber-defense or reliability investment to insurance terms \cite{ref127,ref155}. Evidence from disaster insurance similarly shows that mitigation incentives become stronger when risk-reducing actions can be observed and rewarded through the contract \cite{hudson2016insuranceAdaptation,mol2020moralHazard}.

At a higher level of coupling, insurance terms may be updated according to modeled residual risk or determined jointly with power-system decisions \cite{ref201,ref203}. The key distinction is therefore whether insurance merely transfers realized loss or also responds to actions that reshape future risk. This feedback creates the basis for insurance to influence power-system planning and operation, as examined systematically in Section~\ref{sec:decision-feedback}.

\section{Actuarial assessment and pricing}
\label{sec:actuarial}

\begin{figure*}[width=\textwidth,pos=t]
	\centering
	\includegraphics[width=\textwidth]{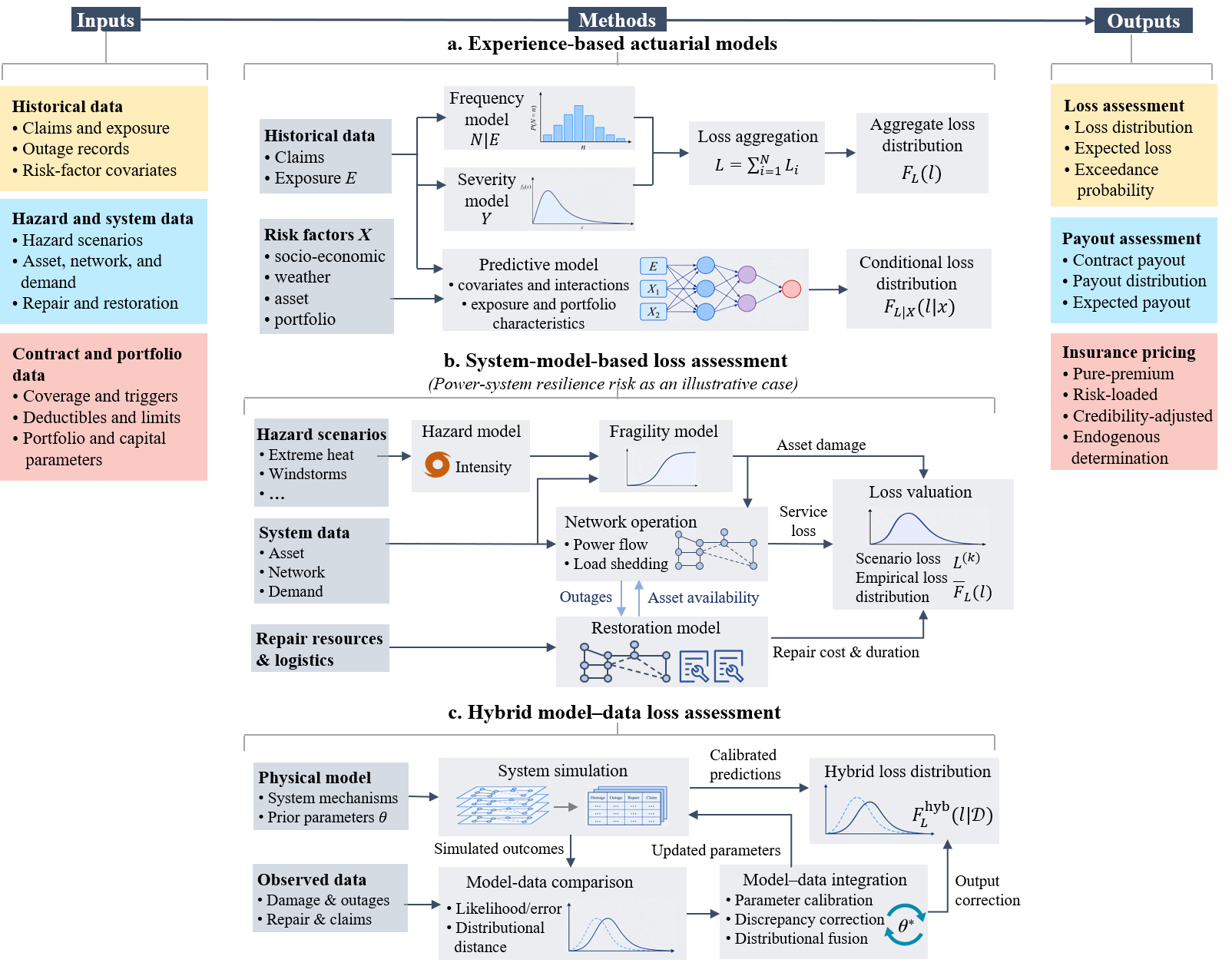}
	\caption{Actuarial assessment and pricing framework for power-system insurance.}
	\label{fig:actuarial-framework}
\end{figure*}

Once the transferred risk has been defined, actuarial analysis focuses on two related tasks: estimating the insured-loss distribution and translating that distribution into an insurance price. Fig.~\ref{fig:actuarial-framework} summarizes this assessment and pricing process.


\subsection{Actuarial approaches to insured-loss assessment}

Let $L$ denote the gross monetary loss associated with the insured interest before contractual terms are applied. The loss distribution and its key summary quantities are represented as
\begin{equation}
F_L(l)=\Pr(L\leq l),
\qquad
\mu_L=\mathbb E[L],
\qquad
\overline F_L(l)=\Pr(L>l),
\label{eq:loss-outputs}
\end{equation}
where $l$ is a loss threshold, $F_L(l)$ is the cumulative distribution function, $\mu_L$ is the expected loss, and $\overline F_L(l)$ is the exceedance probability. The mean characterizes the expected loss level, while the tail of the distribution and dependence across exposures are important for risk loading and capital requirements.  

The following approaches differ primarily in the information used to estimate $L$, ranging from observed loss experience to power-system simulation and combinations of the two.

\subsubsection{Experience-based actuarial models}

Experience-based models infer insured-loss behaviour from observed claims and exposure data. Their main strength is that the modelled response is already expressed in the monetary quantity relevant to insurance. Their reliability depends on whether historical claims are sufficiently representative of the exposure being priced. Two forms are particularly relevant.

\textit{1) Frequency--severity models.}
Claim frequency describes how often covered losses occur, while severity describes the monetary loss conditional on occurrence.  Together, they determine aggregate loss behaviour \cite{awiszus2023cyberPricing}. Frequency should be interpreted relative to exposure, since identical claim counts can imply different levels of risk when the underlying exposure differs. The distinction between frequency and severity is also useful because risk-reduction measures may affect them differently. Preventive actions may primarily reduce the occurrence of loss, whereas operational flexibility and recovery capability may have a stronger influence on loss severity. Historical claims can capture these relationships when sufficiently comparable observations are available, but provide weaker evidence for infrequent or structurally changing risks.

\textit{2) Risk-factor-based models.}
Observed loss distributions can also be conditioned on measurable characteristics of the insured risk,
\begin{equation}
F_{L\mid X}(l\mid x),
\end{equation}
where $X$ denotes observable risk factors, such as weather, asset, and operating characteristics. Generalized linear and hierarchical models provide established structures for relating such information to claim frequency or severity \cite{nelder1997credibility}, while quantile-based methods can capture heterogeneous effects across the loss distribution \cite{eling2022cyberQuantile}.

These models support a priori risk differentiation by using observable characteristics to distinguish expected loss before individual claim experience is incorporated. Credibility adjustment, discussed later, instead updates the rate a posteriori as new experience accumulates. In power-system applications, many predictive models first estimate physical or operational outcomes, which must then be linked to the economic exposure of the insured party. Renewable-insurance studies provide examples of this mapping from resource conditions to project-level financial loss \cite{liao2022solarIndex,mosquera2022weather}, while index-based products additionally require the relationship between the observable trigger and realized loss to be calibrated \cite{niakh2025basis}.

\subsubsection{System-model-based loss assessment}

System-model-based assessment constructs insured-loss realizations through power-system simulation rather than relying primarily on historical claims. For simulation sample $k$,
\begin{equation}
L^{(k)}
=
\mathcal{M}_{L}
\left(
\omega^{(k)},x,u
\right),
\qquad k=1,\ldots,N,
\label{eq:model-based-loss}
\end{equation}
where $\omega^{(k)}$ denotes the uncertain state, $(x,u)$ represents planning and operational decisions, and $\mathcal{M}_{L}$ maps the resulting system consequence to the protected economic loss. Repeated simulation then produces an empirical representation of the insured-loss distribution. The assessment therefore follows a power-system-specific path from uncertain conditions, through system response, to stakeholder-specific economic loss.

The actuarial value of this approach lies in preserving the mechanisms that shape the resulting loss distribution. Correlated initiating outages can materially change cascading-blackout risk \cite{clarfeld2020correlated}, while large-scale weather conditions can simultaneously affect electricity demand and renewable availability \cite{sundar2023meteorological}. Network response further determines how these disturbances translate into service disruption \cite{ref57,ref60,ref67}. The resulting insured loss therefore depends strongly on the system state in which the disturbance occurs and on the response available at that time.

The resulting scenario--loss samples can be summarized through actuarial quantities such as average annual loss and exceedance-probability curves. Average annual loss represents the long-run expected burden, while exceedance-probability curves characterize the upper tail of the loss distribution. Occurrence exceedance probability describes severe individual events, whereas aggregate exceedance probability reflects the accumulation of losses over a specified period. These quantities connect power-system simulation with insurance pricing, reinsurance, and capital assessment \cite{grossi2005catastropheModeling,steptoe2022catModeling}.

Existing applications illustrate this framework across different types of power-system risk. Renewable-energy insurance uses simulated operating and market outcomes to derive producer-specific losses \cite{ref150}, while catastrophe-oriented studies propagate external hazards through system response and recovery before translating the resulting consequences into insured loss \cite{ref202}. The defining feature of system-model-based assessment is the explicit representation of the mechanisms through which uncertain conditions become economic losses.

This approach is particularly useful when observed claims are limited, when historical experience is becoming less representative, or when power-system decisions materially change the loss distribution. Its reliability, however, depends on the credibility of the modelling chain used to generate those losses. Uncertainty in risk representation, system response, recovery, economic valuation, and contract definition can all propagate into the final loss distribution. Model detail is therefore most valuable when it improves the representation of the loss characteristics that matter for pricing and risk-bearing decisions.

\subsubsection{Hybrid model--data loss assessment}

Hybrid assessment combines system-model-generated losses with empirical evidence. This approach is particularly relevant to power systems because physical observations may be relatively abundant even when consistently defined insurance claims remain limited. Observed data can therefore be used to calibrate, correct, or update the loss distribution generated by the system model.

Let $\theta$ denote uncertain model parameters and $\mathcal D$ the available observations. One representative form of parameter-updated hybrid assessment is the Bayesian posterior predictive distribution
\begin{equation}
F_L^{\mathrm{hyb}}(l\mid\mathcal D)
=
\int
F_L(l\mid\theta)\,
p(\theta\mid\mathcal D)
\,\mathrm d\theta ,
\label{eq:hybrid-loss}
\end{equation}
where $F_L(l\mid\theta)$ is the conditional loss distribution for a given parameter realization, $p(\theta\mid\mathcal D)$ is the posterior distribution of the model parameters given the observations, and $F_L^{\mathrm{hyb}}(l\mid\mathcal D)$ is the resulting posterior predictive loss distribution.
When parameter uncertainty alone cannot account for the observed loss behaviour, model discrepancy may be represented explicitly \cite{kennedy2001calibration}.

This review organizes hybrid model--data assessment into three complementary strategies: parameter calibration, model-discrepancy correction, and distributional fusion. Table~\ref{tab:hybrid-strategies} summarizes how each strategy combines simulation and observed evidence in constructing the insured-loss distribution.

\begin{table}[width=\columnwidth,cols=4,pos=t]
	\centering
	\caption{Model--data integration strategies and their actuarial roles in power-system insured-loss assessment.}
	\label{tab:hybrid-strategies}
	
	\fontfamily{ptm}\selectfont
	\scriptsize
	\setlength{\tabcolsep}{3pt}
	\renewcommand{\arraystretch}{1.15}
	
	\begin{tabularx}{\columnwidth}{
			@{}
			>{\hsize=0.85\hsize\linewidth=\hsize
				\RaggedRight\arraybackslash}X
			>{\hsize=1.00\hsize\linewidth=\hsize
				\RaggedRight\arraybackslash}X
			>{\hsize=1.25\hsize\linewidth=\hsize
				\RaggedRight\arraybackslash}X
			>{\hsize=0.90\hsize\linewidth=\hsize
				\RaggedRight\arraybackslash}X
			@{}}
		\toprule
		
		\textbf{Hybrid strategy} &
		\textbf{Model--data interface} &
		\textbf{Actuarial role} &
		\textbf{Main limitation} \\
		\midrule
		
		Parameter calibration
		&
		Parameters governing physical or economic loss formation
		&
		Updates uncertain model parameters before constructing $F_L$
		&
		Structural model error remains
		\\[1.2mm]
		
		Model-discrepancy correction
		&
		Difference between simulated and observed outcomes
		&
		Corrects systematic bias in model-generated losses
		&
		Parameter--discrepancy confounding
		\\[1.2mm]
		
		Distributional fusion
		&
		Observed losses and simulated rare-event losses
		&
		Anchors frequent losses to observations while supplementing the sparsely observed tail with simulations
		&
		Consistency and weighting of heterogeneous evidence
		\\
		
		\bottomrule
	\end{tabularx}
\end{table}

Parameter calibration updates uncertain model parameters before loss scenarios are generated, while discrepancy correction addresses systematic differences between simulated and observed outcomes. Distributional fusion further combines empirical information from more frequently observed losses with simulation-based evidence for sparsely observed tail events. Existing power-system studies already provide relevant methodological foundations, including index recalibration \cite{niakh2025basis}, data-informed estimation of weather-related failure risk \cite{wang2026hierarchicalInsurance}, claims-informed calibration of load-curtailment penalties \cite{liu2026proactiveDefense}, and catastrophe simulation \cite{ref202}. These developments provide useful foundations for a more integrated hybrid assessment of insured-loss distributions.

The three actuarial approaches serve different evidential settings. Experience-based assessment is most effective when comparable claims are sufficiently informative, whereas system-model-based assessment becomes more important when rare or changing risks are poorly represented in historical records. Hybrid approaches connect these two sources of information by using observations to constrain and update model-generated losses. The appropriate approach therefore depends on the evidence available and on the part of the loss distribution that must be represented with confidence.

\subsection{Insurance pricing}

Pricing translates the assessed loss distribution into a monetary charge for risk transfer. This requires assigning a monetary value to uncertain insured payouts. Classical actuarial theory provides a range of valuation principles that account for expected loss and the additional financial implications of uncertainty and tail exposure \cite{goovaerts2003unified,wang2000distortion}. These principles provide the theoretical basis for linking the insured-loss distribution to the premium and the associated risk-bearing requirement.

This review organizes insurance pricing into four connected components. The first three concern the actuarial construction and updating of insurance rates. Pure-premium pricing establishes the expected-payout baseline, risk-loaded pricing incorporates the additional cost of bearing uncertain and concentrated losses, and credibility adjustment updates the rate as relevant insured experience accumulates. Endogenous premium determination then extends this actuarial foundation by allowing insurance terms to respond to power-system decisions and participant behaviour, thereby providing a bridge from actuarial pricing to the decision-coupling mechanisms examined in Section~4.

\subsubsection{Pure-premium pricing}

Let $Z$ denote a contract-relevant observable trigger and $I(L,Z;c)$ the payout under contract configuration $c$. The pure premium is
\begin{equation}
\pi_0(c)
=
\mathbb E\!\left[I(L,Z;c)\right].
\label{eq:pure-premium}
\end{equation}
The pure premium represents the expected contractual payout under the specified coverage structure. Its value therefore depends on how the contract maps the underlying loss and any observable trigger into compensation. Under a frequency--severity formulation, the same principle can be expressed through expected claim frequency and conditional severity relative to the relevant exposure.

Different contracts written on the same underlying physical risk can therefore produce different pure premiums. Contract terms determine the portion of economic loss transferred to the insurer and, for parametric products, the relationship between the observable trigger and realized loss. 
Power-system studies derive the underlying exposure from renewable and market uncertainty \cite{ref150,liao2022solarIndex,mosquera2022weather} as well as from catastrophe and restoration modelling \cite{ref202}. Parametric pricing additionally requires the trigger--loss relationship to be represented with sufficient accuracy to evaluate basis risk.


\subsubsection{Risk-loaded pricing}

Pure premium captures expected contractual payout, while a sustainable insurance rate must also reflect the financial burden arising from loss uncertainty, dependence, and adverse tail outcomes. For an insured portfolio, let \(S\) denote the aggregate contractual payout
\begin{equation}
S
=
\sum_{j=1}^{n}
I(L_j,Z_j;c_j),
\label{eq:portfolio-payout}
\end{equation}
where each term represents the payout associated with insured exposure \(j\). The required risk-bearing capacity is shaped by the tail behaviour of \(S\) and the dependence among its constituent losses. Therefore, portfolios with similar expected payouts can imply very different capital requirements when losses are strongly correlated.

Actuarial theory provides several approaches for assigning monetary value to this distributional risk. Representative formulations are summarized in Table~\ref{tab:risk-valuation}. They differ in how additional weight is assigned to uncertainty or severe loss outcomes \cite{goovaerts2003unified,wang2000distortion,rockafellar2000cvar}. These approaches provide different ways of valuing risk beyond expected payout. Utility-based and Esscher-type principles adjust the monetary value assigned to uncertain or severe losses, while distortion-based methods alter the weighting of the loss distribution. Tail measures such as value-at-risk (VaR), conditional VaR (CVaR), and tail VaR (TVaR) characterize severe portfolio losses that are particularly relevant to capital requirements.

\begin{table}[width=\columnwidth,cols=3,pos=t]
	\centering
	\caption{Representative actuarial risk-valuation approaches relevant to insurance pricing.}
	\label{tab:risk-valuation}
	
	\fontfamily{ptm}\selectfont
	\scriptsize
	\setlength{\tabcolsep}{3pt}
	\renewcommand{\arraystretch}{1.15}
	
	\begin{tabularx}{\columnwidth}{
			@{}
			>{\RaggedRight\arraybackslash}p{1.35cm}
			>{\centering\arraybackslash}p{3.15cm}
			>{\RaggedRight\arraybackslash}X
			@{}}
		\toprule
		
		\textbf{Approach} &
		\textbf{Representative formulation} &
		\textbf{Main pricing implication} \\
		\midrule
		
		Utility-based
		&
		$\begin{aligned}
			&\mathbb{E}[u(W+\pi-S)]\\
			&\qquad =u(W)
		\end{aligned}$
		&
		Reflects risk aversion in the valuation of uncertain losses.
		\\[1.5mm]
		
		Esscher-type
		&
		$\displaystyle
		\pi_{\theta}
		=
		\frac{\mathbb{E}[Se^{\theta S}]}
		{\mathbb{E}[e^{\theta S}]}$
		&
		Places greater weight on large-loss realizations.
		\\[1.5mm]
		
		Distortion-based
		&
		$\displaystyle
		\pi_g
		=
		\int_{0}^{\infty}
		g[\Pr(S>s)]\,\mathrm{d}s$
		&
		Reweights loss probabilities to reflect distributional risk.
		\\[1.5mm]
		
		Tail-risk measures
		&
		$\begin{aligned}
			&\mathrm{TVaR}_{\alpha}(S)\\
			&\quad =
			\frac{1}{1-\alpha}
			\int_{\alpha}^{1}
			\mathrm{VaR}_{q}(S)\,\mathrm{d}q
		\end{aligned}$
		&
		Quantifies severe portfolio losses relevant to capital requirements and risk loading.
		\\
		
		\bottomrule
	\end{tabularx}
	
	\par\vspace{1mm}
	\raggedright
	\fontfamily{ptm}\selectfont
	\scriptsize
	\textit{Note:}
	\(u(\cdot)\) denotes the utility function,
	\(W\) initial wealth,
	\(\pi\) the premium,
	\(\theta\) the Esscher parameter,
	\(g(\cdot)\) a distortion function,
	and \(\alpha\) the confidence level.
	\par
	
\end{table}

Risk loading is therefore inherently portfolio dependent. The capital requirement associated with an individual policy depends on how its losses interact with the insurer’s existing portfolio.  Tail-risk measures can quantify severe portfolio losses \cite{rockafellar2000cvar}, while actuarial pricing  further translates the resulting risk-bearing requirement into  a charge for the transferred risk  \cite{tsanakas2005riskPricing}. Reinsurance and model uncertainty may further affect this cost \cite{kousky2019disasterInsuranceReview,ref175}.

This issue is particularly relevant to power systems, where common disturbances and network-mediated propagation can create strong dependence across insured losses. Similar expected payouts may therefore correspond to very different capital requirements, while diversification across weakly dependent exposures can reduce aggregate risk.

Existing power-system insurance studies have only partially incorporated this portfolio perspective. Premiums are often represented through loading factors, solvency constraints, or predefined risk--premium relationships, with limited connection to the capital implied by correlated loss distributions. A more complete pricing framework would link system-model-generated losses to portfolio dependence, capital requirements, and the resulting cost of risk transfer.

\subsubsection{Credibility-adjusted pricing}

Credibility adjustment provides a formal mechanism for incorporating emerging insured experience into an existing rate,
\begin{equation}
\widehat{\pi}_{t+1}
=
w_t\widehat{\pi}_{\mathrm{exp},t}
+
(1-w_t)\widehat{\pi}_{\mathrm{base},t},
\qquad
0\leq w_t\leq1,
\label{eq:credibility-rate}
\end{equation}
where $\hat{\pi}_{\mathrm{exp},t}$ denotes the premium estimate derived from emerging insured experience, $\hat{\pi}_{\mathrm{base},t}$ denotes the baseline rate based on the broader risk class or prior actuarial assessment, and $w_t$ is the associated credibility weight. The weight generally increases as sufficiently comparable experience accumulates \cite{buhlmann1967credibility}. Hierarchical generalized linear models provide related structures for combining collective experience with individual risk characteristics \cite{nelder1997credibility}. Credibility adjustment therefore updates an existing rate using newly observed insured experience.

Its application to power systems depends critically on the comparability of that experience. Changes in observed claims may reflect changes in the underlying insured risk, but can also arise from shifts in external conditions or system exposure. This distinction becomes more important under non-stationary risk, where historical experience may become less representative of future conditions. Credibility updating should therefore account for changes in both the insured exposure and the environment in which losses are observed.

Power-system studies have begun to incorporate related ideas through updated failure assessment and risk-responsive premium setting \cite{ref127,wang2026hierarchicalInsurance,ref201,ref203}. A fuller credibility-based framework would combine emerging claims with updated information on the underlying risk environment so that rate adjustment remains consistent with the evolving insured exposure.

\subsubsection{Endogenous premium determination}

Endogenous premium determination extends actuarial pricing into a broader decision framework by allowing insurance prices or coverage terms to vary with power-system decisions. In power-system applications, this allows insurance terms to adjust with risk-reduction actions, resource decisions, and participation behaviour \cite{ref127,lau2023mutual,ref203,ref161,xie2023coalitional}. Huang et al. further develop an insurer-optimization framework in which premiums and coverage are jointly determined with end-user participation and the coordination of flexible resources, including vehicle-to-grid resources, to manage compensation exposure \cite{huang2026insurerV2G}.

The resulting premium reflects the interactions represented in the underlying model. Its actuarial interpretation, however, still depends on the insured-loss distribution and the associated risk-bearing requirements used to support the optimized or negotiated terms. Accordingly, endogenous premium determination is most meaningful when these terms are anchored to the expected payout and risk loading established through actuarial assessment.

This feature is particularly relevant to power systems because insurance terms can influence decisions that subsequently reshape the insured-loss distribution. Endogenous pricing therefore creates a feedback between actuarial assessment and power-system decision-making: assessed loss affects the premium and coverage terms, while the resulting decisions alter future risk and can in turn motivate subsequent rate adjustment.

\section{Insurance-aware power-system decision-making}
\label{sec:decision-feedback}

Insurance becomes decision-relevant when power-system actions affect financial exposure or when insurance terms respond to changes in system risk. This section first introduces four levels of insurance--power-system decision coupling and then examines their implications across major decision contexts.

\subsection{Levels of insurance--power-system decision coupling}

Insurance can be incorporated into power-system decision models with different levels of decision coupling. The four coupling levels introduced here characterize the extent to which insurance terms and power-system decisions interact within the decision-making process.

Building on the notation introduced above, let $x$ denote the power-system decision and $c$ the insurance contract configuration. Since $x$ can affect the underlying loss distribution, the resulting system risk may influence both the insurance premium and the loss retained by the system operator. A general insurance-aware decision problem can be formulated as
\begin{equation}
	\begin{aligned}
		\min_{x,c}\quad
		& C(x)
		+\pi\!\left(c;F_{L,Z|x}\right)
		+\rho\!\left[L(x)-I(L(x),Z;c)\right],\\
		\text{s.t.}\quad
		& x\in\mathcal{X}_{\mathrm{security}},
		\qquad
		c\in\mathcal{C}_{\mathrm{contract}} .
	\end{aligned}
	\label{eq:insurance-aware-decision}
\end{equation}
Here, $C(x)$ denotes the power-system decision cost, $\rho(\cdot)$ evaluates the retained loss, and $\mathcal{X}_{\mathrm{security}}$ and $\mathcal{C}_{\mathrm{contract}}$ represent the physical and contractual feasible sets, respectively. The formulation allows the premium to depend on both the contract terms and the loss distribution induced by the power-system decision.

%
%
%

\begin{figure*}[width=\textwidth,pos=t]
	\centering
	\includegraphics[width=0.95\textwidth]{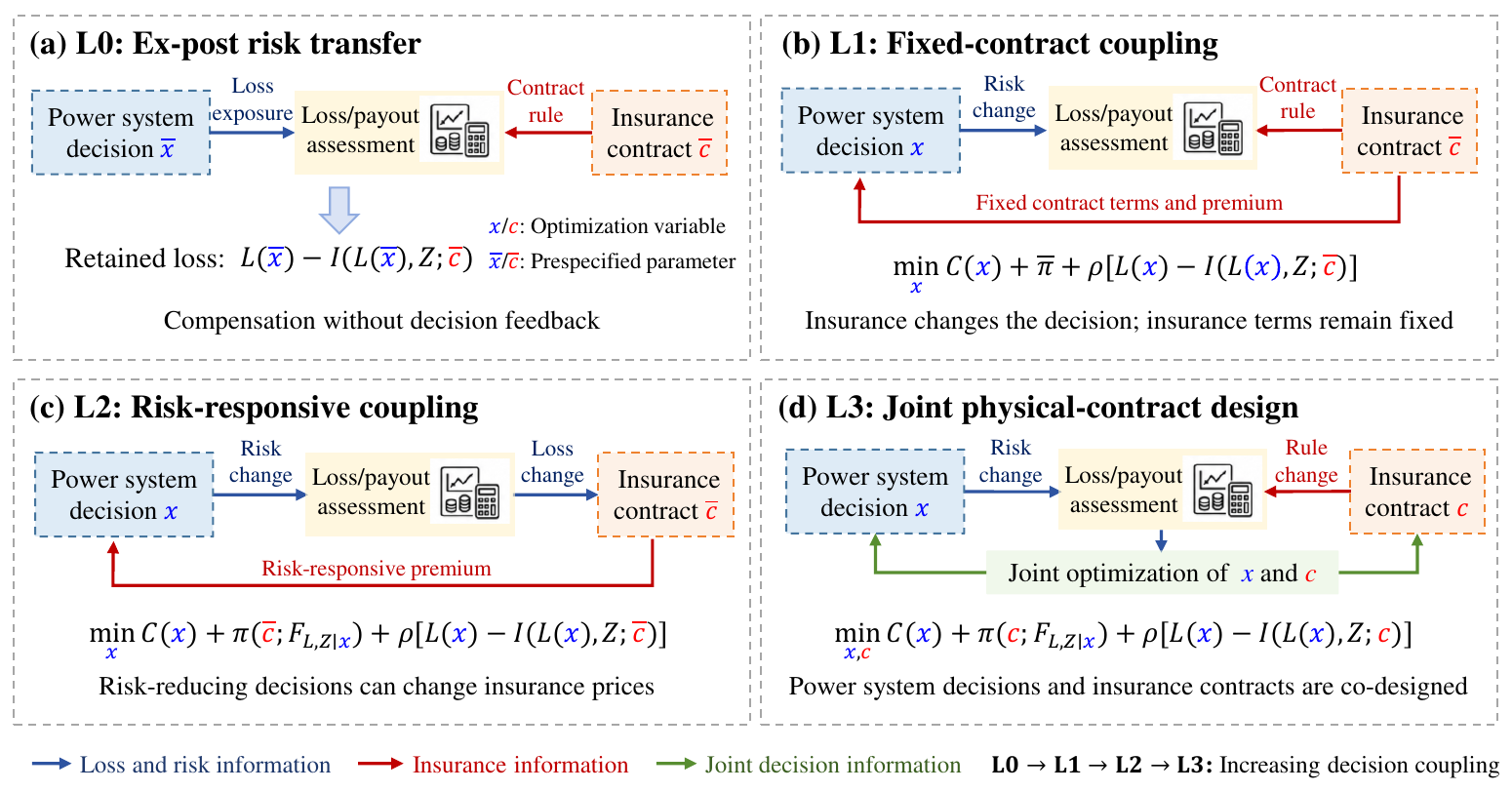}
	\caption{Four levels of insurance--power-system decision coupling, from ex-post risk transfer to joint physical--contract design.}
	\label{fig:coupling}
\end{figure*}

Four levels of insurance--power-system decision coupling are distinguished in Fig.~\ref{fig:coupling}.

\textbf{L0: Ex-post risk transfer.}
The power-system decision is fixed before insurance is evaluated. Realized loss and any relevant trigger are then mapped through the contract to determine the payout. Insurance therefore acts only after the decision has been made, with no feedback to planning or operation.

\textbf{L1: Fixed-contract coupling.}
A predetermined contract $\bar c$ is incorporated into the power-system decision,
\begin{equation}
\min_x\;
C(x)
+\bar{\pi}
+\rho\!\left[
L(x)-I\!\left(L(x),Z;\bar c\right)
\right].
\label{eq:coupling-L1}
\end{equation}
The power-system decision responds to the protection provided by insurance, while the premium and contract terms remain fixed. Insurance can therefore affect investment or operation without responding to the resulting change in risk.

\textbf{L2: Risk-responsive coupling.}
At this level, changes in the power-system decision affect the insured-risk assessment and feed back into insurance terms. A premium-responsive formulation can be written as
\begin{equation}
\min_x\;
C(x)
+\pi\!\left(\bar c;F_{L,Z\mid x}\right)
+\rho\!\left[
L(x)-I\!\left(L(x),Z;\bar c\right)
\right].
\label{eq:coupling-L2}
\end{equation}
The response may involve the premium or other contract conditions. L2 therefore creates a two-way feedback between power-system decisions and financial terms, with the strength of this incentive depending on how credibly changes in risk are reflected in the contract.

\textbf{L3: Joint physical--contract design.}
At the highest level of coupling, both the power-system decision $x$ and contract configuration $c$ are determined jointly according to Eq.~\eqref{eq:insurance-aware-decision}. This allows physical risk reduction and financial transfer to be coordinated within the same decision problem. Recent insurer-oriented models extend this idea by jointly determining insurance provision and the coordination of flexible resources \cite{huang2026insurerV2G}.

\begin{figure*}[width=\textwidth,pos=t]
	\centering
	\includegraphics[width=\textwidth]{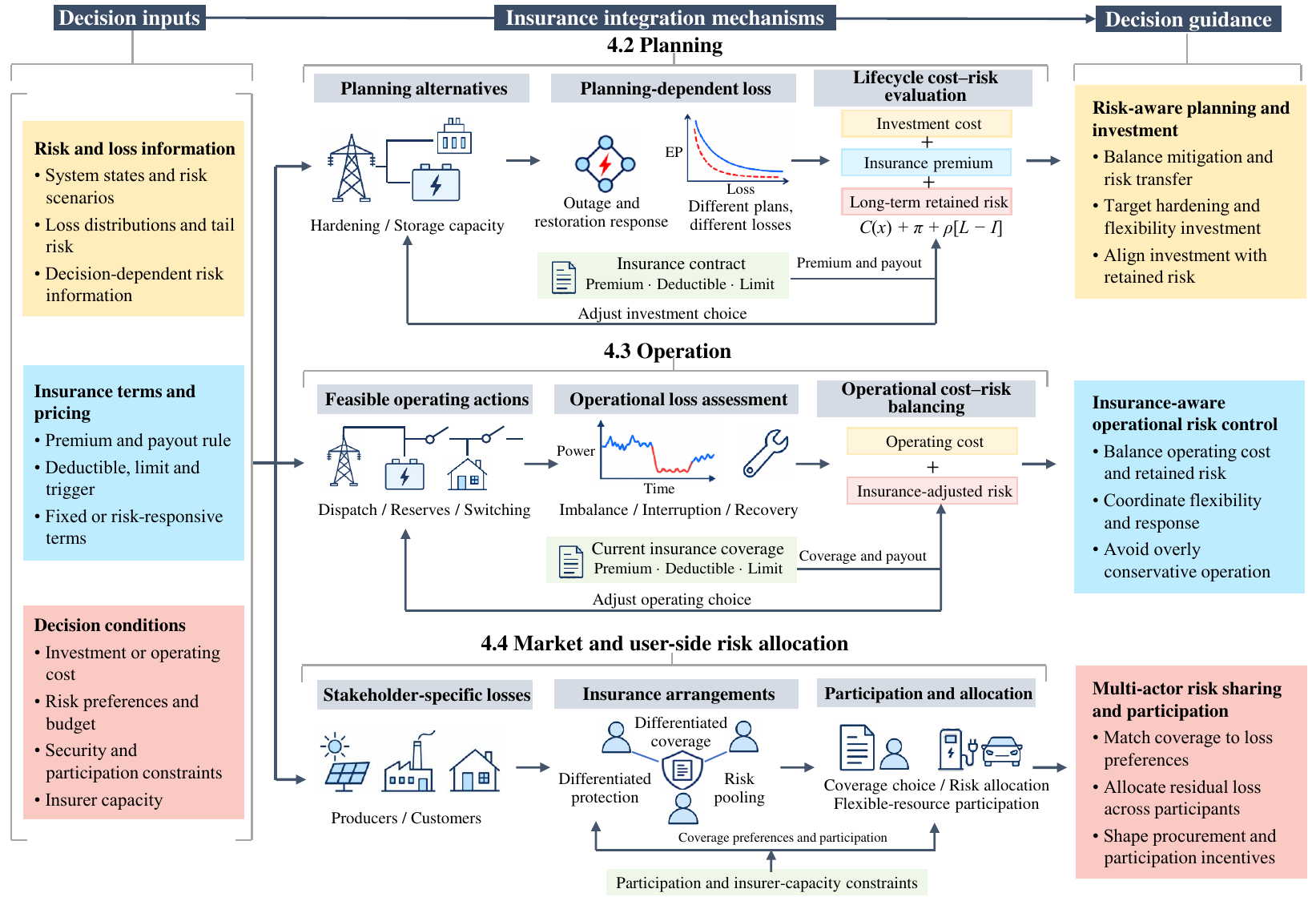}
	\caption{Insurance integration mechanisms in power-system decision-making.}
	\label{fig:decision_framework}
\end{figure*}

Building on this framework, Fig.~\ref{fig:decision_framework} summarizes how insurance is incorporated into major power-system decision contexts, from risk and loss assessment to risk-sensitive decision guidance.


\subsection{Insurance in power-system planning}

Planning decisions shape the long-term distribution of residual loss by changing system exposure and the capacity available to absorb future disturbances. Insurance extends this planning problem by complementing physical risk reduction with financial risk transfer. Their interaction is inherently bidirectional: physical investment can reduce future insured losses and improve insurance conditions, while insurance changes the economic value assigned to additional risk reduction.

Existing studies have begun to embed this interaction directly into planning models. Joint formulations have coordinated physical investment with catastrophe or resilience insurance \cite{ref203,ref201}, while reliability-oriented mechanisms have linked heterogeneous outage exposure to reserve procurement and distributed resilience investment \cite{ref154,ref155}. Similar ideas have been extended to cyber-defense allocation \cite{ref127}, insurance-budget-supported planning \cite{xu2026insuranceBudget}, and insurer-aware infrastructure investment with endogenous pricing and cost allocation \cite{huang2025resilienceInsurancePlanning}. These studies collectively show that planning decisions can reshape insured loss, while insurance terms can in turn alter the value and financing of risk-reduction investment.

A key limitation is the imbalance between physical and financial modelling detail. Power-system planning models increasingly represent network reinforcement, hardening, and flexible resources under uncertain hazards with substantial system detail \cite{ref61,ref76,ref79,ref92,ref100,ghosh2023resilienceInvestment}, while climate-adaptive planning further accounts for long-term changes in resource conditions and climate exposure \cite{gu2026climateAdaptivePlanning}. In contrast, insurance components are typically represented at a much simpler level. The next step is therefore to connect decision-dependent loss distributions with actuarially credible pricing and risk-bearing conditions within the same planning framework.

\subsection{Insurance in operational risk control}

Operational decisions reshape residual loss over short time scales while most infrastructure remains fixed. Insurance changes the economic consequences associated with these decisions and can therefore alter the relative value of feasible operating strategies.

Existing studies have begun to incorporate insurance into operational coordination. For routine uncertainty, insurance has been combined with storage and other flexible resources to manage renewable imbalance exposure and influence scheduling decisions \cite{ref150,ref204,ref152}. Related formulations use insurance contracts to allocate operational responsibilities among independently controlled resources and to influence their participation or dispatch decisions \cite{liu2022dispatchInsurance,ref131}. These studies show how financial protection can be linked directly to short-term resource operation.

Insurance has also been considered under severe disturbances. Recent studies combine risk-informed dispatch with insurance purchasing under weather-related threats \cite{wang2026hierarchicalInsurance}, while insurer-oriented models use insurance financing to coordinate external flexibility during disruptions \cite{huang2026insurerV2G}. Catastrophe-risk models further connect operational and restoration decisions with losses transferred through insurance \cite{ref202}. Together, these studies show how operating decisions reshape insured loss and how insurance terms affect the choice among feasible responses.

A remaining challenge is the mismatch between operational and insurance time scales. Power-system conditions can change rapidly, whereas underwriting and premium adjustment usually occur over much longer intervals. Operational coupling is therefore most naturally built around pre-established coverage and service arrangements, with realized operating experience feeding into subsequent risk assessment and contract updates.

\subsection{Market and user-side risk allocation}

Planning and operation primarily shape the physical consequences that remain after mitigation, whereas market and user arrangements influence how the resulting economic losses are allocated among different parties. Insurance becomes relevant when this allocation changes participation decisions, procurement choices, or the provision of reliability services.

For renewable producers, insurance provides financial protection against production and market uncertainty. Weather- and index-based products provide protection against resource variability \cite{han2019windIndex,liao2022solarIndex}, while more decision-oriented formulations link insurance with market participation and procurement choices \cite{jia2022deviationInsurance,ref161}. Peer-to-peer arrangements further allow residual basis risk to be redistributed across geographically diverse producers \cite{niakh2025basis}. These studies show how insurance can move from passive revenue protection toward influencing renewable-market decisions.

Insurance can also affect the participation of flexible resources. Storage-based contracts have been used to support renewable commitments by reserving flexibility against shortfalls \cite{ref152}, while dispatch insurance links compensation with the operating obligations of distributed resources participating in electricity markets \cite{liu2022dispatchInsurance}. In such settings, insurance changes the value of providing flexibility and can therefore influence resource availability and market behaviour.

On the user side, insurance provides a means of differentiating protection according to heterogeneous interruption exposure and reliability preferences. Existing studies link reliability insurance with reserve procurement and differentiated reliability provision \cite{ref154,ref156}, while peer-to-peer schemes redistribute outage losses among participating users \cite{ref205}. Broader arrangements extend this logic across coupled infrastructures \cite{xie2023coalitional} and among end users, insurers, and external flexibility providers \cite{huang2026insurerV2G}. These studies show how stakeholder-specific loss exposure can be translated into financial incentives for reliability support.

Table~\ref{tab:decision-synthesis} compares representative applications according to their decision focus, transferred exposure, and form of insurance--decision interaction. Across these applications, insurance--decision coupling ranges from externally specified insurance at L1, through risk-responsive interaction at L2, to joint physical--contract design at L3. The appropriate level of coupling depends on the decision context and time scale. Planning can support deeper physical--financial integration, whereas operational and market applications often rely on pre-established contracts with later adjustment. As coupling deepens, credible feedback among power-system decisions, insured losses, and financial terms becomes increasingly important. Existing studies generally model the power-system side in greater detail than the financial side, highlighting insurer feasibility and empirical validation as key areas for further development. At the same time, closer integration makes it necessary to distinguish insurance clearly from related financial and electricity-market mechanisms.

\begin{table*}[width=\textwidth,cols=7,pos=t]
	\centering
	\caption{Critical synthesis of representative insurance-aware power-system decision models.}
	\label{tab:decision-synthesis}
	
	\fontfamily{ptm}\selectfont
	\scriptsize
	\renewcommand{\arraystretch}{1.10}
	\setlength{\tabcolsep}{2.4pt}
	
	\begin{tabularx}{\textwidth}{
			@{}
			>{\RaggedRight\arraybackslash}p{1.45cm}
			>{\RaggedRight\arraybackslash}p{2.25cm}
			>{\RaggedRight\arraybackslash}p{2.05cm}
			>{\RaggedRight\arraybackslash}p{3.15cm}
			>{\RaggedRight\arraybackslash}p{2.05cm}
			>{\centering\arraybackslash}p{0.75cm}
			>{\RaggedRight\arraybackslash}X
			@{}}
		\toprule
		
		\textbf{Decision domain} &
		\textbf{Decision focus} &
		\textbf{Transferred exposure} &
		\textbf{Insurance--decision feedback} &
		\textbf{Representative studies} &
		\textbf{Level} &
		\textbf{Main unresolved boundary} \\
		\midrule
		
		\textbf{Planning}
		&
		Resilience assets and preventive measures
		&
		Catastrophe-induced asset and interruption loss
		&
		Hardening, storage, maintenance, and insurance are jointly selected according to their effects on residual loss.
		&
		\cite{ref201,ref203,huang2025resilienceInsurancePlanning}
		&
		L3
		&
		Premium functions and insurance capacity are generally prescribed rather than derived from portfolio risk.
		\\
		\cmidrule(lr){2-7}
		
		&
		Distributed resilience investment
		&
		Locational and stakeholder-specific outage loss
		&
		Insurance payments or risk-responsive terms convert avoided interruption loss into incentives for local resilience investment.
		&
		\cite{ref155}
		&
		L2--L3
		&
		The mapping from local reliability improvement to premiums, participation, and affordability remains stylized.
		\\
		\cmidrule(lr){2-7}
		
		&
		Insurance-funded and cyber-defense resource allocation
		&
		Extreme-weather and cyber-induced service loss
		&
		Premium-based budgets or cyber-risk charges influence the allocation of generation, storage, reserves, and defense resources.
		&
		\cite{ref127,xu2026insuranceBudget}
		&
		L2--L3
		&
		Insurer funding authority, solvency constraints, and accumulation across correlated exposures remain external.
		\\
		\midrule
		
		\textbf{Operation}
		&
		Renewable-imbalance management
		&
		Forecast-error and imbalance-settlement costs
		&
		Storage absorbs routine deviations, while insurance transfers larger residual losses or supports shared-storage deviation services.
		&
		\cite{ref150,ref204}
		&
		L1--L2
		&
		Coverage and deductibles respond only weakly to the risk reduction achieved through repeated operation.
		\\
		\cmidrule(lr){2-7}
		
		&
		Dispatch compliance and flexible-resource operation
		&
		Non-delivery, dispatch-conflict, and cyber-related operating loss
		&
		Insurance contracts affect storage reservation, distributed-resource compliance, and charging or discharging decisions.
		&
		\cite{ref152,liu2022dispatchInsurance,ref131}
		&
		L1--L2
		&
		Network constraints and operational behaviour are represented more explicitly than underwriting and claim verification.
		\\
		\cmidrule(lr){2-7}
		
		&
		Extreme-event dispatch and restoration
		&
		Weather-induced load loss and emergency-service cost
		&
		Updated failure risk, insurance purchasing, restoration actions, and insurer-funded emergency resources are coordinated before or during disruption.
		&
		\cite{ref202,wang2026hierarchicalInsurance,huang2026insurerV2G}
		&
		L2--L3
		&
		Short-term insurance adjustment is rarely supported by actual underwriting time scales or event-level claims evidence.
		\\
		\midrule
		
		\textbf{Market / user side}
		&
		Renewable production and market participation
		&
		Weather-dependent revenue, production shortfall, and imbalance loss
		&
		Index or deviation coverage stabilizes revenue and can influence insurance purchase, procurement, and market participation.
		&
		\cite{han2019windIndex,liao2022solarIndex,jia2022deviationInsurance,ref161,niakh2025basis}
		&
		L1--L2
		&
		Basis risk, spatial dependence, and the response of market behaviour to premiums remain incompletely validated.
		\\
		\cmidrule(lr){2-7}
		
		&
		Reliability differentiation and local risk pooling
		&
		Consumer-specific interruption loss
		&
		Reliability preferences are translated into differentiated protection, reserve procurement, household insurance choice, or peer-to-peer pooling.
		&
		\cite{ref154,ref156,ref205}
		&
		L2
		&
		Participation, affordability, adverse selection, and the boundary with regulated reliability obligations require stronger treatment.
		\\
		\cmidrule(lr){2-7}
		
		&
		Insurance pricing and multi-party risk allocation
		&
		Infrastructure interruption, revenue loss, and cyber-related loss
		&
		Premiums interact with service prices, while cooperative contracts allocate losses and insurance costs among systems, utilities, or infrastructure operators.
		&
		\cite{xie2023coalitional,ref128,ref130,zhao2026coalitionalTransmission}
		&
		L2--L3
		&
		Pricing feedback and cooperative allocation require stronger treatment of insurer capital, correlated claims, and enforceable settlement rules.
		\\
		
		\bottomrule
	\end{tabularx}
\end{table*}

Table~\ref{tab:boundary} summarizes these boundaries. The key distinction is that insurance transfers a defined insured loss. In contrast, weather derivatives settle against observable weather or production indices \cite{ref157,ref158}, power-purchase agreements mainly stabilize revenue exposure, and capacity mechanisms or reliability options reward availability or scarcity performance \cite{ref165,ref167}. This distinction becomes particularly important when several mechanisms respond to the same scarcity or outage event, as their payment rules and contractual responsibilities should remain clearly defined to avoid overlapping compensation or distorted incentives. Hybrid designs may combine several mechanisms, but the insured interest and payment rules should remain clear. For user-side applications, customer differentiation should also operate within a minimum reliability standard, with insurance addressing residual economic loss rather than substituting for basic network adequacy.

\begin{table*}[width=\textwidth,cols=5,pos=t]
	\centering
	\caption{Boundary between insurance and related financial and electricity-market mechanisms.}
	\label{tab:boundary}
	
	\fontfamily{ptm}\selectfont
	\scriptsize
	\setlength{\tabcolsep}{4pt}
	\renewcommand{\arraystretch}{1.15}
	
	\begin{tabularx}{\textwidth}{
			@{}
			>{\RaggedRight\arraybackslash}p{2.5cm}
			>{\RaggedRight\arraybackslash}p{3.3cm}
			>{\RaggedRight\arraybackslash}p{3.0cm}
			>{\RaggedRight\arraybackslash}p{3.1cm}
			>{\RaggedRight\arraybackslash}X
			@{}}
		\toprule
		
		\textbf{Mechanism} &
		\textbf{Transferred object} &
		\textbf{Settlement basis} &
		\textbf{Power-system role} &
		\textbf{Distinction from insurance} \\
		\midrule
		
		Weather derivative or futures contract
		&
		Tradable weather, price, or production index
		&
		Market-index payoff
		&
		Cash-flow hedging
		&
		Does not require a claimant-specific insured loss
		\\[1.2mm]
		
		Power-purchase agreement
		&
		Energy delivery and price exposure
		&
		Contracted output and price
		&
		Revenue stabilization
		&
		Primarily reallocates commercial price or volume exposure
		\\[1.2mm]
		
		Capacity mechanism
		&
		Qualified resource availability
		&
		Availability and scarcity rules
		&
		Adequacy procurement
		&
		Pays for capability rather than compensating a claimant's loss
		\\[1.2mm]
		
		Reliability option
		&
		Scarcity-price and performance exposure
		&
		Strike-price and performance obligation
		&
		Scarcity-risk control
		&
		Creates a standardized market obligation rather than individualized loss transfer
		\\[1.2mm]
		
		Insurance
		&
		Defined claimant loss or insured index
		&
		Verified loss or agreed trigger
		&
		Residual financial-risk transfer
		&
		Requires a defined coverage relationship and risk-bearing capacity
		\\
		
		\bottomrule
	\end{tabularx}
\end{table*}

\section{Research gaps}
\label{sec:gaps}

Existing studies have established important foundations for transferring residual power-system risk and incorporating insurance into power-system decisions. However, several limitations still restrict the credibility and applicability of current approaches.

\subsection{Limited validation of physical-to-insured-loss translation}


Actuarial assessment requires a credible mapping from a power-system disturbance to the monetary loss covered by insurance. The information needed for this mapping is typically dispersed across physical, operational, market, and insurance records, which are rarely aligned at the level required for event-based validation. This fragmentation makes it difficult to identify where uncertainty enters the physical-to-financial loss chain. Power-system studies increasingly use detailed operating data for outage reconstruction and restoration analysis \cite{ref192,ref193}, but comparable event-level information on insured loss is much less available. As a result, experience-based and system-model-based assessments are seldom evaluated against the same observed loss outcomes. The accuracy of the resulting insured-loss distributions therefore remains difficult to verify, which weakens the actuarial basis for pricing and contract design.

\subsection{Incomplete representation of systemic and non-stationary insurance risk}


Power-system losses are increasingly shaped by dependence across exposures and by changes in the underlying risk environment. Climate evolution and long-term system transformation can change future loss conditions \cite{ref21,ref24,ref27}, while common disturbances and network interactions can generate strongly correlated outcomes \cite{clarfeld2020correlated,ref110}. Modern power-system models increasingly represent these effects, including compound supply--demand stress \cite{turner2019compound,sundar2023meteorological} and network propagation under severe storms \cite{sturmer2024texas}.
The actuarial challenge is to preserve these dependencies as physical consequences are translated into insured losses and aggregated across a portfolio. Correlated system impacts can become simultaneous claims and materially change diversification, capital requirements, and the cost of risk transfer. Existing studies often represent either the physical dependence structure or the financial tail exposure in detail, with a weaker connection between the two. Non-stationarity adds a further difficulty because loss distributions can change with system and environmental conditions. If these effects are not carried consistently into portfolio assessment, pricing and capital requirements may not reflect the actual level and concentration of future risk.

\subsection{Limited endogeneity of insurance in power-system decisions}

Power-system models increasingly account for decisions that alter residual loss, while insurance terms are often introduced through fixed contracts or prescribed pricing rules. More endogenous formulations have begun to connect insurance with resilience investment, maintenance, cyber defense, reliability decisions, and cooperative risk allocation \cite{ref127,ref155,ref201,ref203,xie2023coalitional,huang2026insurerV2G}. These studies establish the basic feedback between risk-reduction decisions and financial transfer.  
The remaining gap lies in the financial side of this interaction. Insurance terms are frequently modeled in less detail than the physical decision process, particularly in pricing, underwriting response, and risk-bearing capacity. As a result, the insurance conditions used in optimization may not be fully supported by the underlying insured-loss distribution or by the insurer's ability to absorb transferred risk.  
Time-scale consistency is another limitation. Long-term planning, contract renewal, and short-term operation evolve on different horizons, but these feedbacks are often treated within a single decision stage or under fixed insurance assumptions. This can overstate how quickly insurance terms respond to operating conditions and weaken the practical interpretation of insurance--power-system coupling.

\section{Conclusions and future perspectives}
\label{sec:future}

This review has examined insurance as a financial complement to physical risk management in power systems, covering insured-risk definition and transfer, actuarial loss assessment and pricing, and insurance--power-system decision coupling. Existing studies show that insurance can not only compensate realized losses but also influence planning, operation, and other risk-reduction decisions through financial incentives. Building on the gaps discussed above, future research should focus on developing more credible, forward-looking, and practical insurance frameworks for power systems. The following subsections outline four key directions.

\subsection{Physical--financial data linkage and loss-model validation}

A priority is to establish event-level evidence that links power-system consequences with the economic losses ultimately recognized by insurance. At present, operational and restoration records are rarely aligned with stakeholder exposure and claims information, which limits direct validation of the physical-to-financial loss chain. Future research should therefore develop data structures that connect these observations across the same events while respecting the proprietary nature of utility and insurance data. Such linkage would allow uncertainty to be traced through the loss-assessment process and would provide a common basis for validating different actuarial approaches. Existing outage and restoration analytics offer a strong physical foundation \cite{ref192,ref193}. Combining them with observed financial outcomes would support more rigorous comparison of experience-based, system-model-based, and hybrid loss models.

\subsection{Systemic and non-stationary portfolio risk}

Future actuarial analysis should represent how power-system dependence propagates into insurance portfolios. Common disturbances can generate correlated losses across many insured exposures, reducing diversification and increasing the capital required to support transferred risk. Portfolio assessment should therefore preserve this dependence from the physical model through insured-loss aggregation and pricing. The underlying loss distribution should also evolve with changing system and environmental conditions \cite{ref21,ref24,ref27}. Forward-looking physical scenarios can be combined with distributionally robust methods \cite{ref180,ref188}, while catastrophe-risk finance provides useful experience in structuring risk-bearing capacity for severe losses \cite{kousky2019disasterInsuranceReview,ref175}. These developments could support more credible assessment of future portfolio tails and the financial capacity needed to absorb them.

\subsection{Multi-timescale insurance--power-system decision coupling}

Future decision models should strengthen the feedback between changes in insured risk and the insurance terms applied to that risk. Risk-responsive pricing and coverage should be supported by observable changes in loss exposure and by a credible representation of underwriting and risk-bearing capacity. The timing of this feedback is equally important. Long-term planning can interact directly with contract design and underwriting, while short-term operation usually takes place under previously established coverage. Operational experience can then inform later reassessment and renewal. Sequential or periodically updated formulations are therefore a promising direction for representing insurance--power-system interaction across different decision horizons.

\subsection{Pilot implementation and empirical evaluation}

Further progress requires applications in which both physical outcomes and insurance performance can be observed. Early pilot studies are likely to be most informative where the insured interest and settlement boundary are clearly defined, making it possible to compare insurance terms with realized risk and subsequent decision responses. Repeated observations across events and renewal cycles would provide evidence on whether insurance remains financially sustainable and whether the resulting incentives persist over time. Such empirical evaluation is essential for identifying which insurance mechanisms can move from conceptual models toward practical power-system applications.

\section*{Declaration of competing interest}
The authors declare that they have no known competing financial interests or personal relationships that could have appeared to influence the work reported in this manuscript.

\section*{Acknowledgements}
This work was supported in part by the National Natural Science Foundation of China under Grant 52607198, in part by the China Postdoctoral Science Foundation under Grant 2025M770483, and in part by the Sichuan Natural Science Foundation Project under Grant                                                                                                            2026NSFSCZY0091.

\bibliographystyle{insurance-unsrt}
\bibliography{references}

@article{awiszus2023cyberPricing,

  title = {Modeling and pricing cyber insurance},
  author  = {Awiszus, Kerstin and Knispel, Thomas and Penner, Irina and Svindland, Gregor and Vo{\ss}, Alexander and Weber, Stefan},
  journal = {European Actuarial Journal},
  volume  = {13},
  number  = {1},
  pages   = {1--53},
  year    = {2023},
}

@article{goovaerts2003unified,

  title = {A unified approach to generate risk measures},
  author = {Goovaerts, Marc J and Kaas, Rob and Dhaene, Jan and Tang, Qihe},
  journal={ASTIN Bulletin: The Journal of the IAA},
  volume={33},
  number={2},
  pages={173--191},
  year={2003},
}

@article{wang2000distortion,

  title = {A class of distortion operators for pricing financial and insurance risks},
  author = {Wang, Shaun S.},
  journal = {Journal of Risk and Insurance},
  volume  = {67},
  number  = {1},
  pages   = {15--36},
  year    = {2000},
}

@article{rockafellar2000cvar,

  title = {Optimization of conditional value-at-risk},
  author = {Rockafellar, R. Tyrrell and Uryasev, Stanislav},
  journal = {Journal of Risk},
  volume  = {2},
  number  = {3},
  pages   = {21--41},
  year    = {2000},
}

@article{ghosh2023resilienceInvestment,

  title = {A stochastic investment decision making method for distribution system resilience enhancement considering automation, hardening and distributed energy resources},
  author = {Ghosh, Puspendu and De, Mala},
  journal = {Reliability Engineering \& System Safety},
  volume  = {237},
  pages   = {109395},
  year    = {2023},
}

@article{han2019windIndex,

  title = {Weather index insurance for wind energy},
  author = {Han, Xiao and Zhang, Guoxing and Xie, Yixiang and Yin, Jiaxuan and Zhou, Haiming and Yang, Yunxi and Li, Junhui and Bai, Wenlei},
  journal = {Global Energy Interconnection},
  volume  = {2},
  number  = {6},
  pages   = {541--548},
  year    = {2019},
}

@article{liao2022solarIndex,

  title = {Index-based renewable energy insurance for {Taiwan} solar photovoltaic power plants},
  author = {Liao, S.-C. and Chang, S.-C. and Cheng, T.-C.},
  journal = {Risk Management and Insurance Review},
  volume  = {25},
  number  = {2},
  pages   = {145--172},
  year    = {2022},
}

@article{clement2018basisRisk,

  title = {A global review of the impact of basis risk on the functioning of and demand for index insurance},
  author = {Clement, Kristina Y. and Botzen, W. J. Wouter and Brouwer, Roy and Aerts, Jeroen C. J. H.},
  journal = {International Journal of Disaster Risk Reduction},
  volume  = {28},
  pages   = {845--853},
  year    = {2018},
}

@article{gorman2022vollReview,

  title = {The quest to quantify the value of lost load: a critical review of the economics of power outages},
  author = {Gorman, Will},
  journal = {The Electricity Journal},
  volume  = {35},
  number  = {8},
  pages   = {107187},
  year    = {2022},
}

@article{dionne2022catCapacity,

  title = {A re-examination of the {US} insurance market's capacity to pay catastrophe losses},
  author = {Dionne, Georges and Desjardins, Denise},
  journal = {Risk Management and Insurance Review},
  volume  = {25},
  number  = {4},
  pages   = {515--549},
  year    = {2022},
}

@article{lakdawalla2012catBondReinsurance,

  title = {Catastrophe bonds, reinsurance, and the optimal collateralization of risk transfer},
  author = {Lakdawalla, Darius and Zanjani, George},
  journal = {Journal of Risk and Insurance},
  volume  = {79},
  number  = {2},
  pages   = {449--476},
  year    = {2012},
}

@article{clarke2017sovereignRiskFinance,

  title = {{Evaluating sovereign disaster risk finance strategies: A framework}},
  author = {Clarke, Daniel J. and Mahul, Olivier and Poulter, Richard and Teh, Tse-Ling},
  journal = {The Geneva Papers on Risk and Insurance - Issues and Practice},
  volume  = {42},
  pages   = {565--584},
  year    = {2017},
}

@article{zhang2024blended,

  title ={{Blended insurance scheme: A synergistic conventional-index insurance mixture}},
  author = {Zhang, Jinggong},
  journal = {Insurance: Mathematics and Economics},
  volume  = {119},
  pages   = {93--105},
  year    = {2024},
}

@article{lopez2026indexHybrid,

  author = {Lopez, Olivier and Nkameni, Daniel},

  title = {Index insurance under demand and solvency constraints},
  journal = {European Actuarial Journal},
  year    = {2026},
  volume  = {16},
  pages   = {297--328},
}

@article{poudyal2026resiliencePlanning,

  title = {Resilience-driven planning of electric power systems against extreme weather events},
  author = {Poudyal, Abodh and Lamichhane, Shishir and do Prado, Josue Campos and Dubey, Anamika},
  journal = {IET Generation, Transmission \& Distribution},
  volume  = {20},
  number  = {1},
  pages   = {e70330},
  year    = {2026},
}

@article{ciapessoni2023costBenefit,

  title = {A cost--benefit analysis framework for power system resilience enhancement based on optimization via simulation considering climate changes and cascading outages},
  author = {Ciapessoni, Emanuele and Cirio, Diego and Pitto, Andrea},
  journal = {Energies},
  volume  = {16},
  number  = {13},
  pages   = {5160},
  year    = {2023},
}

@book{berliner1982limits,

  author = {Berliner, Baruch},
  title     = {Limits of Insurability of Risks},
  publisher = {Prentice-Hall},
  address   = {Englewood Cliffs, NJ},
  year      = {1982},
}

@article{berliner1985large,

  author = {Berliner, Baruch},

  title = {Large risks and limits of insurability},
  journal = {The Geneva Papers on Risk and Insurance},
  year    = {1985},
  volume  = {10},
  number  = {37},
  pages   = {313--329},
}

@article{biener2015cyber,

  author = {Biener, Christian and Eling, Martin and Wirfs, Jan Hendrik},

  title ={{Insurability of cyber risk: An empirical analysis}},
  journal = {The Geneva Papers on Risk and Insurance - Issues and Practice},
  year    = {2015},
  volume  = {40},
  number  = {1},
  pages   = {131--158},
}

@article{charpentier2008insurability,

  author = {Charpentier, Arthur},

  title = {Insurability of climate risks},
  journal = {The Geneva Papers on Risk and Insurance - Issues and Practice},
  year    = {2008},
  volume  = {33},
  number  = {1},
  pages   = {91--109},
}

@article{vazquez2022wildfire,

  author = {Zuniga Vazquez, Daniel A. and Qiu, Feng and Fan, Neng and Sharp, Kevin},

  title = {{Wildfire mitigation plans in power systems: A literature review}},
  journal = {IEEE Transactions on Power Systems},
  year    = {2022},
  volume  = {37},
  number  = {5},
  pages   = {3540--3551},
}

@article{liao2022outage,

  author = {Liao, Yizheng and Weng, Yang and Tan, Chin-Woo and Rajagopal, Ram},

  title = {Quick line outage identification in urban distribution grids via smart meters},
  journal = {CSEE Journal of Power and Energy Systems},
  year    = {2022},
  volume  = {8},
  number  = {4},
  pages   = {1074--1086},
}

@article{ovaere2019voll,

  author = {Ovaere, Marten and Heylen, Evelyn and Proost, Stef and Deconinck, Geert and Van Hertem, Dirk},

  title = {How detailed value of lost load data impact power system reliability decisions},
  journal = {Energy Policy},
  year    = {2019},
  volume  = {132},
  pages   = {1064--1075},
}

@article{staffell2018weather,

  author = {Staffell, Iain and Pfenninger, Stefan},

  title = {The increasing impact of weather on electricity supply and demand},
  journal = {Energy},
  year    = {2018},
  volume  = {145},
  pages   = {65--78},
}

@article{sundar2023meteorological,

  author = {Sundar, Srihari and Craig, Michael T. and Payne, Ashley E. and Brayshaw, David J. and Lehner, Flavio},

  title = {Meteorological drivers of resource adequacy failures in current and high renewable western {U.S.} power systems},
  journal = {Nature Communications},
  year    = {2023},
  volume  = {14},
  pages   = {6379},
}

@article{clarfeld2020correlated,

  author = {Clarfeld, Laurence A. and Hines, Paul D. H. and Hernandez, Eric M. and Eppstein, Margaret J.},

  title = {Risk of cascading blackouts given correlated component outages},
  journal = {IEEE Transactions on Network Science and Engineering},
  year    = {2020},
  volume  = {7},
  number  = {3},
  pages   = {1133--1144},
}

@article{yang2020premium,

  author = {Yang, Zhiyuan and Liu, Yun and Campbell, Meghan and Ten, Chee-Wooi and Rho, Yeonwoo and Wang, Lingfeng and Wei, Wei},

  title = {Premium calculation for insurance businesses based on cyber risks in {IP}-based power substations},
  journal = {IEEE Access},
  year    = {2020},
  volume  = {8},
  pages   = {78890--78900},
}

@article{lau2023mutual,

  author = {Lau, Pikkin and Wang, Lingfeng and Wei, Wei and Liu, Zhaoxi and Ten, Chee-Wooi},

  title = {A novel mutual insurance model for hedging against cyber risks in power systems deploying smart technologies},
  journal = {IEEE Transactions on Power Systems},
  year    = {2023},
  volume  = {38},
  number  = {1},
  pages   = {630--642},
}

@article{wang2022resilienceReview,

  author = {Wang, Chong and Ju, Ping and Wu, Feng and Pan, Xueping and Wang, Zhaoyu},

  title = {A systematic review on power system resilience from the perspective of generation, network, and load},
  journal = {Renewable and Sustainable Energy Reviews},
  volume  = {167},
  pages   = {112567},
  year    = {2022},
}

@article{younesi2022resilienceReview,

  author = {Younesi, Abdollah and Shayeghi, Hossein and Wang, Zongjie and Siano, Pierluigi and Mehrizi-Sani, Ali and Safari, Amin},

  title = {Trends in modern power systems resilience: state-of-the-art review},
  journal = {Renewable and Sustainable Energy Reviews},
  volume  = {162},
  pages   = {112397},
  year    = {2022},
}

@article{shi2022distributionResilienceReview,

  author = {Shi, Qingxin and Liu, Wenxia and Zeng, Bo and Hui, Hongxun and Li, Fangxing},

  title = {Enhancing distribution system resilience against extreme weather events: concept review, algorithm summary, and future vision},
  journal = {International Journal of Electrical Power \& Energy Systems},
  volume  = {138},
  pages   = {107860},
  year    = {2022},
}

@article{panteli2015weatherResilience,

  author = {Panteli, Mathaios and Mancarella, Pierluigi},

  title = {Influence of extreme weather and climate change on the resilience of power systems: impacts and possible mitigation strategies},
  journal = {Electric Power Systems Research},
  volume  = {127},
  pages   = {259--270},
  year    = {2015},
}

@article{cronin2018climateEnergyReview,

  author = {Cronin, Jennifer and Anandarajah, Gabrial and Dessens, Olivier},

  title = {Climate change impacts on the energy system: a review of trends and gaps},
  journal = {Climatic Change},
  volume  = {151},
  number  = {2},
  pages   = {79--93},
  year    = {2018},
}

@article{bromberger2026extremeWeatherReview,

  author = {Bromberger, Dani Antonini and Garcia, Vinicius Jacques and Rodrigues, St{\'e}fane Dias and Neves, Luciane Silva},

  title = {A taxonomy for extreme weather events and their impacts on power systems: a comprehensive review of classification, effects, and mitigation strategies},
  journal = {Renewable and Sustainable Energy Reviews},
  volume  = {227},
  pages   = {116498},
  year    = {2026},
}

@article{ghanbarzadeh2025capacityPlanningReview,

  author = {Ghanbarzadeh, Taraneh and Habibi, Daryoush and Aziz, Asma},

  title = {Addressing reliability challenges in generation capacity planning under high penetration of renewable energy resources and storage solutions: a review},
  journal = {Renewable and Sustainable Energy Reviews},
  volume  = {212},
  pages   = {115461},
  year    = {2025},
}

@article{haugen2023uncertaintyReview,

  author = {Haugen, Mari and Farahmand, Hossein and Jaehnert, Stefan and Fleten, Stein-Erik},

  title = {Representation of uncertainty in market models for operational planning and forecasting in renewable power systems: a review},
  journal = {Energy Systems},
  year    = {2023},
}

@article{wang2025uncertaintyReview,

  author = {Wang, Zhaoyuan and Bu, Siqi and Wen, Jiaxin and Huang, Can},

  title = {A comprehensive review on uncertainty modeling methods in modern power systems},
  journal = {International Journal of Electrical Power \& Energy Systems},
  volume  = {166},
  pages   = {110534},
  year    = {2025},
}

@article{kousky2019disasterInsuranceReview,

  author = {Kousky, Carolyn},

  title = {The role of natural disaster insurance in recovery and risk reduction},
  journal = {Annual Review of Resource Economics},
  volume  = {11},
  pages   = {399--418},
  year    = {2019},
}

@article{ingels2024climateInsuranceReview,

  author = {Ingels, Michiel W. and Botzen, W. J. Wouter and Aerts, Jeroen C. J. H. and Brusselaers, Jan and Tesselaar, Max},

  title = {The state of the art and future of climate risk insurance modeling},
  journal = {Annals of the New York Academy of Sciences},
  volume  = {1541},
  number  = {1},
  pages   = {100--114},
  year    = {2024},
}

@article{zhao2021catastropheReview,

  author = {Zhao, Yang and Lee, Jin-Ping and Yu, Min-Teh},

  title = {Catastrophe risk, reinsurance and securitized risk-transfer solutions: a review},
  journal = {China Finance Review International},
  volume  = {11},
  number  = {4},
  pages   = {449--473},
  year    = {2021},
}

@article{ref42,

  author = {Kozarcanin, S. and Liu, H. and Andresen, G. B.},

  title = {21st century climate change impacts on key properties of a large-scale renewable-based electricity system},
  journal = {Joule},
  year    = {2019},
  volume  = {3},
  number  = {4},
  pages   = {992--1005},
}

@article{ref21,

  author = {Panteli, M. and Pickering, C. and Wilkinson, S. and Dawson, R. and Mancarella, P.},

  title = {Power system resilience to extreme weather: fragility modeling, probabilistic impact assessment, and adaptation measures},
  journal = {IEEE Transactions on Power Systems},
  year    = {2017},
  volume  = {32},
  number  = {5},
  pages   = {3747--3757},
}

@article{ref23,

  author = {Perera, A. T. D. and Nik, V. M. and Chen, D. and Scartezzini, J. L. and Hong, T.},

  title = {Quantifying the impacts of climate change and extreme climate events on energy systems},
  journal = {Nature Energy},
  year    = {2020},
  volume  = {5},
  pages   = {150--159},
}

@article{ref121,

  author = {Ashok, A. and Govindarasu, M. and Wang, J.},

  title = {Cyber-physical attack-resilient wide-area monitoring, protection, and control for the power grid},
  journal = {Proceedings of the IEEE},
  year    = {2017},
  volume  = {105},
  number  = {7},
  pages   = {1389--1407},
}

@article{ref127,

  author = {Lau, P. and Wei, W. and Wang, L. and Liu, Z. and Ten, C.-W.},

  title = {A cybersecurity insurance model for power system reliability considering optimal defense resource allocation},
  journal = {IEEE Transactions on Smart Grid},
  year    = {2020},
  volume  = {11},
  number  = {5},
  pages   = {4403--4414},
}

@article{ref154,

  author = {Billimoria, F. and Fele, F. and Savelli, I. and Morstyn, T. and McCulloch, M.},

  title = {An insurance mechanism for electricity reliability differentiation under deep decarbonization},
  journal = {Applied Energy},
  year    = {2022},
  volume  = {321},
  pages   = {119356},
}

@article{ref155,

  author = {Billimoria, F. and Fele, F. and Savelli, I. and Morstyn, T. and McCulloch, M.},

  title = {An insurance paradigm for improving power system resilience via distributed investment},
  journal = {IEEE Transactions on Energy Markets, Policy and Regulation},
  year    = {2023},
  volume  = {1},
  number  = {4},
  pages   = {499--511},
}

@article{ref4,

  author = {Jufri, F. H. and Widiputra, V. and Jung, J.},

  title = {State-of-the-art review on power grid resilience to extreme weather events: definitions, frameworks, quantitative assessment methodologies, and enhancement strategies},
  journal = {Applied Energy},
  year    = {2019},
  volume  = {239},
  pages   = {1049--1065},
}

@article{ref7,

  author = {Umunnakwe, A. and Huang, H. and Oikonomou, K. and Davis, K. R.},

  title = {Quantitative analysis of power systems resilience: standardization, categorizations, and challenges},
  journal = {Renewable and Sustainable Energy Reviews},
  year    = {2021},
  volume  = {149},
  pages   = {111252},
}

@article{ref16,

  author = {Ghiasi, M. and Niknam, T. and Wang, Z. and Mehrandezh, M. and Dehghani, M. and Ghadimi, N.},

  title ={{A comprehensive review of cyber-attacks and defense mechanisms for improving security in smart grid energy systems: Past, present and future}},
  journal = {Electric Power Systems Research},
  year    = {2023},
  volume  = {215},
  pages   = {108975},
}

@article{ref123,

  author = {Liang, G. and Zhao, J. and Luo, F. and Weller, S. R. and Dong, Z. Y.},

  title = {A review of false data injection attacks against modern power systems},
  journal = {IEEE Transactions on Smart Grid},
  year    = {2017},
  volume  = {8},
  number  = {4},
  pages   = {1630--1638},
}

@article{grams2017spatialWind,

  author = {Grams, Christian M. and Beerli, Remo and Pfenninger, Stefan
             and Staffell, Iain and Wernli, Heini},

  title = {Balancing {Europe's} wind-power output through spatial deployment informed by weather regimes},
  journal = {Nature Climate Change},
  year    = {2017},
  volume  = {7},
  pages   = {557--562},
}

@article{gu2026climateAdaptivePlanning,

  author = {Gu, Chenjia and Ruan, Jiaqi and Qiu, Yiwei and Zang, Tianlei
             and Chen, Shi and Xu, Zhao and Wen, Fushuan and Zhang, Pei
             and Dong, Zhao Yang and Wang, Peng},

  title = {{Toward climate-adaptive low-carbon power system planning: A multistage stochastic framework considering climate uncertainties}},
  journal = {IEEE Transactions on Industrial Informatics},
  year    = {2026},
  volume  = {22},
  number  = {5},
  pages   = {3938--3949},
}

@article{gu2026climateAdaptiveReview,

  author = {Gu, Chenjia and Ruan, Jiaqi and Yang, Xinyi and Huang, Jing
             and Qiu, Yiwei and Wang, Jianxiao and Chen, Shi and Liang, Gaoqi
             and Xu, Zhao and Su, Hui and Lie, Tek Tjing and Song, Jie},

  title = {Advancing climate-adaptive operation and planning for renewable-rich energy systems},
  journal = {Renewable and Sustainable Energy Reviews},
  year    = {2026},
  volume  = {242},
  pages   = {117336},
}

@article{ref146,

  author = {Woods, Daniel W. and Moore, Tyler},

  title = {Does insurance have a future in governing cybersecurity?},
  journal = {IEEE Security \& Privacy},
  year    = {2020},
  volume  = {18},
  number  = {1},
  pages   = {21--27},
}

@article{ref147,

  author = {Marotta, A. and Martinelli, F. and Nanni, S. and Orlando, A. and Yautsiukhin, A.},

  title = {Cyber-insurance survey},
  journal = {Computer Science Review},
  year    = {2017},
  volume  = {24},
  pages   = {35--61},
}

@article{ref205,

  author = {Wang, H. and Yang, J. and Li, G. and Wu, F. and Lin, Y. and Liao, F. and Jia, L.},

  title = {Outage risk management for residential users considering distribution network reliability assessment: {A P2P} outage insurance approach},
  journal = {Electric Power Systems Research},
  year    = {2025},
  volume  = {245},
  pages   = {111590},
}

@article{ref2,

  author = {Panteli, M. and Trakas, D. N. and Mancarella, P. and Hatziargyriou, N. D.},

  title = {Power systems resilience assessment: hardening and smart operational enhancement strategies},
  journal = {Proceedings of the IEEE},
  year    = {2017},
  volume  = {105},
  number  = {7},
  pages   = {1202--1213},
}

@article{ref3,

  author = {Bie, Z. and Lin, Y. and Li, G. and Li, F.},

  title = {Battling the extreme: a study on the power system resilience},
  journal = {Proceedings of the IEEE},
  year    = {2017},
  volume  = {105},
  number  = {7},
  pages   = {1253--1266},
}

@article{ref8,

  author = {Mishra, D. K. and Ghadi, M. J. and Azizivahed, A. and Li, L. and Zhang, J.},

  title = {A review on resilience studies in active distribution systems},
  journal = {Renewable and Sustainable Energy Reviews},
  year    = {2021},
  volume  = {135},
  pages   = {110201},
}

@article{ref10,

  author = {Stankovi{\'c}, Aleksandar M. and Tomsovic, Kevin L. and De Caro, Fabrizio and Braun, Martin and Chow, Joe H. and {\v C}ukalevski, Ninel and Dobson, Ian and Eto, Joseph and Fink, Blair and Hachmann, Christian and Hill, David and Ji, Chuanyi and Kavicky, James A. and Levi, Victor and Liu, Chen-Ching and Mili, Lamine and Moreno, Rodrigo and Panteli, Mathaios and Petit, Frederic D. and Sansavini, Giovanni and Singh, Chanan and Srivastava, Anurag K. and Strunz, Kai and Sun, Hongbo and Xu, Yin and Zhao, Shijia},

  title = {Methods for analysis and quantification of power system resilience},
  journal = {IEEE Transactions on Power Systems},
  year    = {2023},
  volume  = {38},
  number  = {5},
  pages   = {4774--4787},
}

@article{ref15,

  author = {Shuva Paul and Abodh Poudyal and Shiva Poudel and Anamika Dubey and Zhaoyu Wang},

  title = {Resilience assessment and planning in power distribution systems: past and future considerations},
  journal = {Renewable and Sustainable Energy Reviews},
  year    = {2024},
  volume  = {189},
  pages   = {113991},
}

@article{ref18,

  author = {Gholami, A. and Shekari, T. and Amirioun, M. H. and Aminifar, F. and Amini, M. H. and Sargolzaei, A.},

  title = {Toward a consensus on the definition and taxonomy of power system resilience},
  journal = {IEEE Access},
  year    = {2018},
  volume  = {6},
  pages   = {32035--32053},
}

@article{ref22,

  author = {Panteli, M. and Mancarella, P. and Trakas, D. N. and Kyriakides, E. and Hatziargyriou, N. D.},

  title = {Metrics and quantification of operational and infrastructure resilience in power systems},
  journal = {IEEE Transactions on Power Systems},
  year    = {2017},
  volume  = {32},
  number  = {6},
  pages   = {4732--4742},
}

@article{ref150,

  author = {Yang, H. and Qiu, J. and Meng, K. and Zhao, J. H. and Dong, Z. Y. and Lai, M.},

  title = {Insurance strategy for mitigating power system operational risk introduced by wind power forecasting uncertainty},
  journal = {Renewable Energy},
  year    = {2016},
  volume  = {89},
  pages   = {606--615},
}

@article{ref151,

  author = {D'Achiardi, David and Aguiar, Nayara and Baros, Stefanos and Gupta, Vijay and Annaswamy, Anuradha M.},

  title = {Reliability contracts between renewable and natural gas power producers},
  journal = {IEEE Transactions on Control of Network Systems},
  year    = {2019},
  volume  = {6},
  number  = {3},
  pages   = {1075--1085},
}

@article{ref29,

  author = {Solaun, K. and Cerd{\'a}, E.},

  title = {Climate change impacts on renewable energy generation. {A} review of quantitative projections},
  journal = {Renewable and Sustainable Energy Reviews},
  year    = {2019},
  volume  = {116},
  pages   = {109415},
}

@article{ref30,

  author = {Tobin, I. and Greuell, W. and Jerez, S. and Ludwig, F. and Vautard, R. and van Vliet, M T H and Br{\'e}on, F-M},

  title = {Vulnerabilities and resilience of {European} power generation to 1.5\textdegree{}{C}, 2\textdegree{}{C} and 3\textdegree{}{C} warming},
  journal = {Environmental Research Letters},
  year    = {2018},
  volume  = {13},
number  = {4},
  pages   = {044024},
}

@article{ref32,

  author = {Burillo, D. and Chester, M. V. and Ruddell, B. and Johnson, N.},

  title = {Electricity demand planning forecasts should consider climate non-stationarity to maintain reserve margins during heat waves},
  journal = {Applied Energy},
  year    = {2017},
  volume  = {206},
  pages   = {267--277},
}

@article{ref38,

  author = {van Ruijven, B. J. and De Cian, E. and Sue Wing, Ian},

  title = {Amplification of future energy demand growth due to climate change},
  journal = {Nature Communications},
  year    = {2019},
  volume  = {10},
  pages   = {2762},
}

@article{ref45,

  author = {Raynaud, D. and Hingray, B. and Fran{\c c}ois, B. and Creutin, J. D.},

  title = {Energy droughts from variable renewable energy sources in {European} climates},
  journal = {Renewable Energy},
  year    = {2018},
  volume  = {125},
  pages   = {578--589},
}

@article{ref46,

  author = {van der Most, Lieke and van der Wiel, Karin and Benders, R. M. J. and Gerbens-Leenes, P. W. and Bintanja, Richard},

  title = {Temporally compounding energy droughts in {European} electricity systems with hydropower},
  journal = {Nature Energy},
  year    = {2024},
  volume  = {9},
  pages={1474--1484},
}

@article{ref25,

  author = {Miara, A. and Macknick, J. E. and V{\"o}r{\"o}smarty, C. J. and Tidwell, V. C. and Newmark, R. and Fekete, B.},

  title = {Climate and water resource change impacts and adaptation potential for {US} power supply},
  journal = {Nature Climate Change},
  year    = {2017},
  volume  = {7},
  pages   = {793--798},
}

@article{ref27,

  author = {Michael T. Craig and Stuart Cohen and Jordan Macknick and Caroline Draxl and Omar J. Guerra and Manajit Sengupta and Sue Ellen Haupt and Bri-Mathias Hodge and Carlo Brancucci},

  title = {A review of the potential impacts of climate change on bulk power system planning and operations in the {United States}},
  journal = {Renewable and Sustainable Energy Reviews},
  year    = {2018},
  volume  = {98},
  pages   = {255--267},
}

@article{ref52,

  author = {Hui Hou and Zhiwei Zhang and Ruizeng Wei and Yong Huang and Yongchao Liang and Xianqiang Li},

  title = {Review of failure risk and outage prediction in power system under wind hazards},
  journal = {Electric Power Systems Research},
  year    = {2022},
  volume  = {210},
  pages   = {108098},
}

@article{ref54,

  author = {Rahman, M. S. and Wahl, T. and Nagaraj, M. and Enr{\'i}quez, A. R.},

  title = {Characteristics of power outages from compound weather extremes in {Florida}},
  journal = {Environmental Research: Climate},
  year    = {2025},
  volume  = {4},
  number  = {3},
  pages   = {035008},
}

@article{ref202,

  author = {Sun, S. and Li, G. and Bian, Y. and Bie, Z. and Hu, Q.},

  title = {Catastrophe risk management for electric power distribution systems: an insurance approach},
  journal = {CSEE Journal of Power and Energy Systems},
  year    = {2023},
  volume  = {9},
  number  = {1},
  pages   = {393--410},
}

@article{zhao2026coalitionalTransmission,

  author = {Zhao, Yuhong and Hu, Zhengyang and Ding, Yibo and Shi, Wenzhuo and Zhang, Aoxiang and Li, Xudong and Xu, Zhao and Chen, Chen and Bie, Zhaohong},

  title = {A coalitional insurance framework for risk management of interconnected transmission systems against extreme weather events},
  journal = {IEEE Transactions on Industrial Informatics},
  year    = {2026},
  volume  = {22},
  number  = {8},
  pages   = {7254--7265},
}

@article{ref130,

  author = {Acharya, S. and Mieth, R. and Konstantinou, C. and Karri, R. and Dvorkin, Y.},

  title = {Cyber insurance against cyberattacks on electric vehicle charging stations},
  journal = {IEEE Transactions on Smart Grid},
  year    = {2022},
  volume  = {13},
  number  = {2},
  pages   = {1529--1541},
}

@article{ref122,

  author = {Liu, X. and Li, Z.},

  title = {False data attacks against {AC} state estimation with incomplete network information},
  journal = {IEEE Transactions on Smart Grid},
  year    = {2017},
  volume  = {8},
  number  = {5},
  pages   = {2239--2248},
}

@inproceedings{ref125,

  author = {Tan, Rui and Nguyen, Hoang Hai and Foo, Eddy Y. S. and Dong, Xinshu and Yau, David K. Y. and Kalbarczyk, Zbigniew  and Iyer, Ravishankar K. and Gooi, Hoay Beng},

  title = {Optimal false data injection attack against automatic generation control in power grids},
  booktitle = {2016 ACM/IEEE 7th International Conference on Cyber-Physical Systems (ICCPS)},
  year      = {2016},
  pages = {1--10},
}

@article{ref82,

  author = {Wang, Yifei and Huang, Liping and Shahidehpour, Mohammad and Lai, Loi Lei and Zhou, Ya},

  title = {Impact of cascading and common-cause outages on resilience-constrained optimal economic operation of power systems},
  journal = {IEEE Transactions on Smart Grid},
  year    = {2020},
  volume  = {11},
  number  = {1},
  pages   = {590--601},
}

@article{ref65,

  author = {Shao, C. and Shahidehpour, M. and Wang, X. and Wang, X. and Wang, B.},

  title = {Integrated planning of electricity and natural gas transportation systems for enhancing the power grid resilience},
  journal = {IEEE Transactions on Power Systems},
  year    = {2017},
  volume  = {32},
  number  = {6},
  pages   = {4418--4429},
}

@article{ref110,

  author = {Fang, Y.-P. and Zio, E.},

  title = {An adaptive robust framework for the optimization of the resilience of interdependent infrastructures under natural hazards},
  journal = {European Journal of Operational Research},
  year    = {2019},
  volume  = {276},
  number  = {3},
  pages   = {1119--1136},
}

@article{ref100,

  author = {Yang, H. and Yang, H. and Rhodes, N. and Roald, L. and Ntaimo, L.},

  title = {Multistage stochastic program for mitigating power system risks under wildfire disruptions},
  journal = {Electric Power Systems Research},
  year    = {2024},
  volume  = {234},
  pages   = {110773},
}

@article{ref117,

  author = {Khalili, Tohid and Davoudi, Masoud and Bidram, Ali},

  title = {Minimizing the reliability trade-off in wildfires risk mitigation},
  journal = {Applied Energy},
  year    = {2024},
  volume  = {371},
  pages   = {123623},
}

@article{ref200,

  author = {Zio, E.},

  title = {Challenges in the vulnerability and risk analysis of critical infrastructures},
  journal = {Reliability Engineering \& System Safety},
  year    = {2016},
  volume  = {152},
  pages   = {137--150},
}

@article{ref198,

  author = {Cremer, Frank and Sheehan, Barry and Fortmann, Michael and Kia, Arash N. and Mullins, Martin and Murphy, Finbarr and Materne, Stefan},

  title = {Cyber risk and cybersecurity: a systematic review of data availability},
  journal = {The Geneva Papers on Risk and Insurance - Issues and Practice},
  year    = {2022},
  volume  = {47},
  number  = {3},
  pages   = {698--736},
}

@article{ref34,

  author = {Bartos, M. and Chester, M. and Johnson, N. and Gorman, B. and Eisenberg, D. and Linkov, I. and Bates, M.},

  title = {Impacts of rising air temperatures on electric transmission ampacity and peak electricity load in the {United States}},
  journal = {Environmental Research Letters},
  year    = {2016},
  volume  = {11},
 number  = {11},
  pages   = {114008},
}

@article{ref24,

  author = {van Vliet, M. T. H. and Wiberg, D. and Leduc, S. and Riahi, K.},

  title = {Power-generation system vulnerability and adaptation to changes in climate and water resources},
  journal = {Nature Climate Change},
  year    = {2016},
  volume  = {6},
  pages   = {375--380},
}

@article{ref37,

  author = {Gernaat, D. E. H. J. and de Boer, H. S. and Daioglou, V. and Yalew, S. G. and M{\"u}ller, C. and van Vuuren, D. P.},

  title = {Climate change impacts on renewable energy supply},
  journal = {Nature Climate Change},
  year    = {2021},
  volume  = {11},
  pages   = {119--125},
}

@article{ref40,

  author = {Ravestein, P. and van der Schrier, G. and Haarsma, R. and Scheele, R. and van den Broek, M.},

  title = {Vulnerability of {European} intermittent renewable energy supply to climate change and climate variability},
  journal = {Renewable and Sustainable Energy Reviews},
  year    = {2018},
  volume  = {97},
  pages   = {497--508},
}

@article{ref174,

  author = {Hudson, P. and Botzen, W. J. W. and Czajkowski, J. and Kreibich, H.},

  title = {Moral hazard in natural disaster insurance markets: empirical evidence from {Germany} and the {United States}},
  journal = {Land Economics},
  year    = {2017},
  volume  = {93},
  number  = {2},
  pages   = {179--208},
}

@article{ref57,

  author = {Chen, C. and Wang, J. and Qiu, F. and Zhao, D.},

  title = {Resilient distribution system by microgrids formation after natural disasters},
  journal = {IEEE Transactions on Smart Grid},
  year    = {2016},
  volume  = {7},
  number  = {2},
  pages   = {958--966},
}

@article{ref67,

  author = {Lei, S. and Wang, J. and Chen, C. and Hou, Y.},

  title = {Mobile emergency generator pre-positioning and real-time allocation for resilient response to natural disasters},
  journal = {IEEE Transactions on Smart Grid},
  year    = {2018},
  volume  = {9},
  number  = {3},
  pages   = {2030--2041},
}

@article{ref75,

  author = {Wang, Ying and Xu, Yin and He, Jinghan and Liu, Chen-Ching and Schneider, Kevin P. and Hong, Mingguo and Ton, Dan T.},

  title = {Coordinating multiple sources for service restoration to enhance resilience of distribution systems},
  journal = {IEEE Transactions on Smart Grid},
  year    = {2019},
  volume  = {10},
  number  = {5},
  pages   = {5781--5793},
}

@article{ref109,

  author = {Yi-Ping Fang and Giovanni Sansavini},

  title = {Optimum post-disruption restoration under uncertainty for enhancing critical infrastructure resilience},
  journal = {Reliability Engineering \& System Safety},
  year    = {2019},
  volume  = {185},
  pages   = {1--11},
}

@article{ref111,

  author = {Espinoza, S. and Panteli, M. and Mancarella, P. and Rudnick, H.},

  title = {Multi-phase assessment and adaptation of power systems resilience to natural hazards},
  journal = {Electric Power Systems Research},
  year    = {2016},
  volume  = {136},
  pages   = {352--361},
}

@article{ref175,

  author = {Braun, A.},

  title ={{Pricing in the primary market for cat bonds: New empirical evidence}},
  journal = {Journal of Risk and Insurance},
  year    = {2016},
  volume  = {83},
  number  = {4},
  pages   = {811--847},
}

@article{ref157,

  author = {Kanamura, T. and Homann, L. and Prokopczuk, M.},

  title = {Pricing analysis of wind power derivatives for renewable energy risk management},
  journal = {Applied Energy},
  year    = {2021},
  volume  = {304},
  pages   = {117827},
}

@article{ref152,

  author = {Aguiar, N. and Gupta, V.},

  title = {An insurance contract design to boost storage participation in the electricity market},
  journal = {IEEE Transactions on Sustainable Energy},
  year    = {2021},
  volume  = {12},
  number  = {1},
  pages   = {543--552},
}

@article{ref165,

  author = {Bublitz, A. and Keles, D. and Zimmermann, F. and Fraunholz, C. and Fichtner, W.},

  title = {A survey on electricity market design: insights from theory and real-world implementations of capacity remuneration mechanisms},
  journal = {Energy Economics},
  year    = {2019},
  volume  = {80},
  pages   = {1059--1078},
}

@article{ref167,

  author = {Paolo Mastropietro and Pablo Rodilla and Michel Rivier and Carlos Batlle},

  title = {{Reliability options: Regulatory recommendations for the next generation of capacity remuneration mechanisms}},
  journal = {Energy Policy},
  year    = {2024},
  volume  = {185},
  pages   = {113959},
}

@article{ref128,

  author = {Lau, P. and Wang, L. and Liu, Z. and Wei, W. and Ten, C.-W.},

  title = {A coalitional cyber-insurance design considering power system reliability and cyber vulnerability},
  journal = {IEEE Transactions on Power Systems},
  year    = {2021},
  volume  = {36},
  number  = {6},
  pages   = {5512--5524},
}

@article{ref203,

  author = {Hu, Q. and Li, G. and Sun, S. and Bie, Z.},

  title = {Incorporating catastrophe insurance in power distribution systems investment and planning for resilience enhancement},
  journal = {International Journal of Electrical Power \& Energy Systems},
  year    = {2024},
  volume  = {155},
  pages   = {109438},
}

@article{ref201,

  author = {Tsao, Y. C. and Thanh, V. V.},

  title = {Integrated voltage control and maintenance insurance planning for distribution networks considering uncertainties},
  journal = {Electric Power Systems Research},
  year    = {2021},
  volume  = {201},
  pages   = {107501},
}

@article{ref204,

  author = {Xu, X. and Li, G. and Yang, H. and Liu, D. and Wang, J. and Zhang, Z. and Zhao, L. and Zhang, X.},

  title = {Pricing method of shared energy storage bias insurance service based on large number theorem},
  journal = {Journal of Energy Storage},
  year    = {2023},
  volume  = {69},
  pages   = {107726},
}

@article{ref184,

  author = {Vrakopoulou, Maria and Margellos, Kostas and Lygeros, John and Andersson, Goran},

  title = {A probabilistic framework for reserve scheduling and {N-1} security assessment of systems with high wind power penetration},
  journal = {IEEE Transactions on Power Systems},
  year    = {2013},
  volume  = {28},
  number  = {4},
  pages   = {3885--3896},
}

@article{ref185,

  author = {Duan, Chao and Jiang, Lin and Fang, Wanliang and Liu, Jun and Liu, Shiming},

  title = {Data-driven distributionally robust energy-reserve-storage dispatch},
  journal = {IEEE Transactions on Industrial Informatics},
  year    = {2018},
  volume  = {14},
  number  = {7},
  pages   = {2826--2836},
}

@article{ref51,

  author = {Dehghani, F. and Mohammadi, M. and Karimi, M.},

  title = {Age-dependent resilience assessment and quantification of distribution systems under extreme weather events},
  journal = {International Journal of Electrical Power \& Energy Systems},
  year    = {2023},
  volume  = {150},
  pages   = {109089},
}

@article{ref139,

  author = {Deng, Ruilong and Xiao, Gaoxi and Lu, Rongxing and Liang, Hao and Vasilakos, Athanasios V.},

  title = {{False data injection on state estimation in power systems--Attacks, impacts, and defense: A survey}},
  journal = {IEEE Transactions on Industrial Informatics},
  year    = {2017},
  volume  = {13},
  number  = {2},
  pages   = {411--423},
}

@article{ref124,

  author = {Sridhar, S. and Govindarasu, M.},

  title = {Model-based attack detection and mitigation for automatic generation control},
  journal = {IEEE Transactions on Smart Grid},
  year    = {2014},
  volume  = {5},
  number  = {2},
  pages   = {580--591},
}

@article{ref60,

  author = {Panteli, M. and Trakas, D. N. and Mancarella, P. and Hatziargyriou, N. D.},

  title = {Boosting the power grid resilience to extreme weather events using defensive islanding},
  journal = {IEEE Transactions on Smart Grid},
  year    = {2016},
  volume  = {7},
  number  = {6},
  pages   = {2913--2922},
}

@article{ref61,

  author = {Yuan, W. and Wang, J. and Qiu, F. and Chen, C. and Kang, C. and Zeng, B.},

  title = {Robust optimization-based resilient distribution network planning against natural disasters},
  journal = {IEEE Transactions on Smart Grid},
  year    = {2016},
  volume  = {7},
  number  = {6},
  pages   = {2817--2826},
}

@article{ref73,

  author = {Ma, Shanshan and Li, Shiyang and Wang, Zhaoyu and Qiu, Feng},

  title = {Resilience-oriented design of distribution systems},
  journal = {IEEE Transactions on Power Systems},
  year    = {2019},
  volume  = {34},
  number  = {4},
  pages   = {2880--2891},
}

@article{ref58,

  author = {Farzin, H. and Fotuhi-Firuzabad, M. and Moeini-Aghtaie, M.},

  title = {Enhancing power system resilience through hierarchical outage management in multi-microgrids},
  journal = {IEEE Transactions on Smart Grid},
  year    = {2016},
  volume  = {7},
  number  = {6},
  pages   = {2869--2879},
}

@article{ref85,

  author = {Hussain, A. and Bui, V-H and Kim, H-M},

  title = {Microgrids as a resilience resource and strategies used by microgrids for enhancing resilience},
  journal = {Applied Energy},
  year    = {2019},
  volume  = {240},
  pages   = {56--72},
}

@article{ref86,

  author = {Galvan, Eric and Mandal, Paras and Sang, Yuanrui},

  title = {Networked microgrids with roof-top solar {PV} and battery energy storage to improve distribution grids resilience to natural disasters},
  journal = {International Journal of Electrical Power \& Energy Systems},
  year    = {2020},
  volume  = {123},
  pages   = {106239},
}

@article{ref156,

  author = {Fuentes, R. and Sengupta, A.},

  title = {{Using insurance to manage reliability in the distributed electricity sector: Insights from an agent-based model}},
  journal = {Energy Policy},
  year    = {2020},
  volume  = {139},
  pages   = {111251},
}

@article{ref88,

  author = {Xu, Jianchun and Zhang, Tingting and Du, Yaxin and Zhang, Wen and Yang, Tianyou and Qiu, Jifu},

  title = {Islanding and dynamic reconfiguration for resilience enhancement of active distribution systems},
  journal = {Electric Power Systems Research},
  year    = {2020},
  volume  = {189},
  pages   = {106749},
}

@article{ref131,

  author = {Hoang, Dinh Thai and Wang, Ping and Niyato, Dusit and Hossain, Ekram},

  title = {Charging and discharging of plug-in electric vehicles ({PEVs}) in vehicle-to-grid ({V2G}) systems: {A} cyber insurance-based model},
  journal = {IEEE Access},
  year    = {2017},
  volume  = {5},
  pages   = {732--754},
}

@article{ref17,

  author = {Li, Y. and Yan, J.},

  title = {{Cybersecurity of smart inverters in the smart grid: A survey}},
  journal = {IEEE Transactions on Power Electronics},
  year    = {2023},
  volume  = {38},
  number  = {2},
  pages   = {2364--2383},
}

@article{ref76,

  author = {Lagos, Tomas and Moreno, Rodrigo and Espinosa, Alejandro Navarro and Panteli, Mathaios and Sacaan, Rafael and Ordonez, Fernando and Rudnick, Hugh and Mancarella, Pierluigi},

  title = {Identifying optimal portfolios of resilient network investments against natural hazards, with applications to earthquakes},
  journal = {IEEE Transactions on Power Systems},
  year    = {2020},
  volume  = {35},
  number  = {2},
  pages   = {1411--1421},
}

@article{ref79,

  author = {Nazemi, Mostafa and Moeini-Aghtaie, Moein and Fotuhi-Firuzabad, Mahmud and Dehghanian, Payman},

  title = {Energy storage planning for enhanced resilience of power distribution networks against earthquakes},
  journal = {IEEE Transactions on Sustainable Energy},
  year    = {2020},
  volume  = {11},
  number  = {2},
  pages   = {795--806},
}

@article{ref92,

  author = {Hou, Hui and Tang, Junyi and Zhang, Zhiwei and Wang, Zhuo and Wei, Ruizeng and Wang, Lei and He, Huan and Wu, Xixiu},

  title = {Resilience enhancement of distribution network under typhoon disaster based on two-stage stochastic programming},
  journal = {Applied Energy},
  year    = {2023},
  volume  = {338},
  pages   = {120892},
}

@article{ref107,

  author = {Zare-Bahramabadi, Majid and Abbaspour, Ali and Fotuhi-Firuzabad, Mahmud and Moeini-Aghtaie, Moein},

  title = {Resilience-based framework for switch placement problem in power distribution systems},
  journal = {IET Generation, Transmission \& Distribution},
  year    = {2018},
  volume  = {12},
  number  = {5},
  pages   = {1223--1230},
}

@article{ref192,

  author = {Wang, Dingwei and Yuan, Yuxuan and Cheng, Rui and Wang, Zhaoyu},

  title = {Data-driven outage restoration time prediction via transfer learning with cluster ensembles},
  journal = {IEEE Transactions on Power Systems},
  year    = {2024},
  volume  = {39},
  number  = {1},
  pages   = {83--96},
}

@article{ref193,

  author = {Ji, C. and Wei, Y. and Poor, H. V.},

  title = {{Resilience of energy infrastructure and services: Modeling, data analytics, and metrics}},
  journal = {Proceedings of the IEEE},
  year    = {2017},
  volume  = {105},
  number  = {7},
  pages   = {1354--1366},
}

@article{ref197,

  author = {Kwasinski, A.},

  title = {Quantitative model and metrics of electrical grids' resilience evaluated at a power distribution level},
  journal = {Energies},
  year    = {2016},
  volume  = {9},
number  = {2},
  pages   = {93},
}

@article{ref80,

  author = {Sabouhi, H. and Doroudi, A. and Fotuhi-Firuzabad, M. and Bashiri, M.},

  title = {Electrical power system resilience assessment: a comprehensive approach},
  journal = {IEEE Systems Journal},
  year    = {2020},
  volume  = {14},
  number  = {2},
  pages   = {2643--2652},
}

@article{ref81,

  author = {Poudel, Shiva and Dubey, Anamika and Bose, Anjan},

  title = {Risk-based probabilistic quantification of power distribution system operational resilience},
  journal = {IEEE Systems Journal},
  year    = {2020},
  volume  = {14},
  number  = {3},
  pages   = {3506--3517},
}

@article{ref180,

  author = {Mohajerin Esfahani, P. and Kuhn, D.},

  title = {{Data-driven distributionally robust optimization using the Wasserstein metric: Performance guarantees and tractable reformulations}},
  journal = {Mathematical Programming},
  year    = {2018},
  volume  = {171},
  pages   = {115--166},
}

@article{ref188,

  author = {Poolla, Bala Kameshwar and Hota, Ashish R. and Bolognani, Saverio and Callaway, Duncan S. and Cherukuri, Ashish},

  title = {{Wasserstein} distributionally robust look-ahead economic dispatch},
  journal = {IEEE Transactions on Power Systems},
  year    = {2021},
  volume  = {36},
  number  = {3},
  pages   = {2010--2022},
}

@article{ref161,

  author = {Xin, Baogui and Hao, Yaru and Peng, Wei and Tan, Hui and Sun, Minghe and Santibanez Gonzalez, Ernesto D. R.},

  title = {Can index-based weather insurance enhance the sustainability of renewable electricity systems? {A} game-theoretic analysis},
  journal = {Reliability Engineering \& System Safety},
  year    = {2026},
  volume  = {272},
  number  = {Part 2},
  pages   = {112577},
}

@article{jia2022deviationInsurance,

  author = {Jia, Qiangang and Chen, Sijie and Li, Yiyan and Yan, Zheng},

  title = {Deviation insurance for risk-averse wind power producers in the {Nordic} power market},
  journal = {International Journal of Electrical Power \& Energy Systems},
  year    = {2022},
  volume  = {134},
  pages   = {107431},
}

@article{zhao2026multiStakeholder,

  author = {Zhao, Yuhong and Ding, Yibo and Sun, Xianzhuo and Shi, Wenzhuo and Hu, Zhengyang and Xu, Zhao and Chen, Chen and Bie, Zhaohong},

  title = {A multi-stakeholder actuarial risk management framework for transmission systems against extreme weather events},
  journal = {IEEE Transactions on Energy Markets, Policy and Regulation},
  year    = {2026},
  volume  = {4},
  number  = {1},
  pages   = {78--92},
}

@article{amonkar2025compositeIndex,

  author = {Amonkar, Yash and Pahel-Short, Corey and Zeighami, Amir and Su, Yufei and Kern, Jordan D. and Characklis, Gregory W.},

  title = {A composite index-based insurance instrument for managing the financial risk of variable hydrometeorology for electric utilities},
  journal = {Environmental Research: Energy},
  year    = {2025},
  volume  = {2},
  number  = {3},
  pages   = {035006},
}

@article{nematshahi2025wildfireCatBond,

  author = {Nematshahi, Saeed and Sohrabi, Behrouz and Arabnya, Ali and Khodaei, Amin and Belval, Erin},

  title = {A catastrophe bond design for the financial resilience of electric utilities against wildfires},
  journal = {IEEE Transactions on Energy Markets, Policy and Regulation},
  year    = {2025},
  volume  = {3},
  number  = {1},
  pages   = {133--143},
}

@article{ref114,

  author = {Peyman Afzali and Farshid Keynia and Masoud Rashidinejad},

  title = {A new model for reliability-centered maintenance prioritisation of distribution feeders},
  journal = {Energy},
  year    = {2019},
  volume  = {171},
  pages   = {701--709},
}

@inproceedings{huang2025resilienceInsurancePlanning,

  title = {Resilience-oriented planning and cost allocation of energy storage integrated with soft open point based on resilience insurance},
  author = {Huang, Bingkai and Huang, Yuxiong and Hu, Qianwen and Li, Gengfeng and Bie, Zhaohong},
  booktitle={2025 IEEE Power \& Energy Society General Meeting (PESGM)},
  pages={1--5},
  year={2025},
  organization={IEEE}
}

@article{liu2022dispatchInsurance,

  author = {Liu, Yikui and Wu, Lei and Shahidehpour, Mohammad},

  title = {Dispatch insurance for {ISO}-market participating {DERs} in distribution systems},
  journal = {IEEE Transactions on Power Systems},
  year    = {2022},
  volume  = {37},
  number  = {4},
  pages   = {3165--3168},
}

@article{liu2026proactiveDefense,

  author = {Liu, Jiangtao and Gao, Jianwei and Meng, Qichen
             and Tao, Hongyuan and Wu, Haoyu},

  title = {{Proactive defense of multi-energy microgrids under compound meteorological hazards: An insurance-actuarial-data-driven {GLUEVaR} framework}},
  journal = {Applied Energy},
  year    = {2026},
  volume  = {422},
  pages   = {128285},
}

@article{ref115,

  author = {Brian Jaramillo-Leon and Jonatas B. Leite},

  title = {Multi-objective optimization for preventive tree trimming scheduling in overhead electric power distribution networks},
  journal = {Journal of Control, Automation and Electrical Systems},
  year    = {2022},
  volume  = {33},
  pages   = {550--560},
}

@article{ref135,

  author = {Liu, Mengxiang and Teng, Fei and Zhang, Zhenyong and Ge, Pudong and Sun, Mingyang and Deng, Ruilong and Cheng, Peng and Chen, Jiming},

  title = {Enhancing cyber-resiliency of {DER}-based smart grid: {A} survey},
  journal = {IEEE Transactions on Smart Grid},
  year    = {2024},
  volume  = {15},
  number  = {5},
  pages   = {4998--5030},
}

@article{ref94,

  author = {Hao Wu and Yunyun Xie and Yan Xu and Qiuwei Wu and Chen Yu and Jinsheng Sun},

  title = {Resilient scheduling of {MESSs} and {RCs} for distribution system restoration considering the forced cut-off of wind power},
  journal = {Energy},
  year    = {2022},
  volume  = {244},
  pages   = {123081},
}

@article{ref112,

  author = {Diego Piasson and Andr{\'e} A. P. B{\'i}scaro and F{\'a}bio B. Le{\~a}o and Jos{\'e} Roberto Sanches Mantovani},

  title = {A new approach for reliability-centered maintenance programs in electric power distribution systems based on a multiobjective genetic algorithm},
  journal = {Electric Power Systems Research},
  year    = {2016},
  volume  = {137},
  pages   = {41--50},
}

@article{ref158,

  author = {Gersema, G. and Wozabal, D.},

  title = {An equilibrium pricing model for wind power futures},
  journal = {Energy Economics},
  year    = {2017},
  volume  = {65},
  pages   = {64--74},
}

@article{turner2019compound,

  author = {Turner, S. W. D. and Voisin, N. and Fazio, J. and Hua, D. and Jourabchi, M.},

  title = {Compound climate events transform electrical power shortfall risk in the {Pacific Northwest}},
  journal = {Nature Communications},
  year    = {2019},
  volume  = {10},
  pages   = {8},
}

@article{sturmer2024texas,

  author = {St{\"u}rmer, Julian and Plietzsch, Anton and Vogt, Thomas and Hellmann, Frank and Kurths, J{\"u}rgen and Otto, Christian and Frieler, Katja and Anvari, Mehrnaz},

  title = {Increasing the resilience of the {Texas} power grid against extreme storms by hardening critical lines},
  journal = {Nature Energy},
  year    = {2024},
  volume  = {9},
  pages   = {526--535},
}

@article{xie2023coalitional,

  author = {Xie, Haipeng and Sun, Xiaotian and Fu, Wei and Chen, Chen and Bie, Zhaohong},

  title = {Risk management for integrated power and natural gas systems against extreme weather: a coalitional insurance contract approach},
  journal = {Energy},
  year    = {2023},
  volume  = {263},
  pages   = {125750},
}

@article{mosquera2022weather,

  author = {Mosquera-L{\'o}pez, Stephania and Uribe, Jorge M.},

  title = {Pricing the risk due to weather conditions in small variable renewable energy projects},
  journal = {Applied Energy},
  year    = {2022},
  volume  = {322},
  pages   = {119476},
}

@article{niakh2025basis,

  author = {Niakh, Fallou and Bassi{\`e}re, Alicia and Denuit, Michel and Robert, Christian Y.},

  title = {Peer-to-peer basis risk management for renewable production parametric insurance},
  journal = {Annals of Operations Research},
  year    = {2025},
}

@article{buhlmann1967credibility,

  author = {B{\"u}hlmann, Hans},

  title = {Experience rating and credibility},
  journal = {ASTIN Bulletin: The Journal of the IAA},
  year    = {1967},
  volume  = {4},
  number  = {3},
  pages   = {199--207},
}

@article{nelder1997credibility,

  author = {Nelder, J. A. and Verrall, R. J.},

  title = {Credibility theory and generalized linear models},
  journal = {ASTIN Bulletin: The Journal of the IAA},
  year    = {1997},
  volume  = {27},
  number  = {1},
  pages   = {71--82},
}

@article{kennedy2001calibration,

  author = {Kennedy, Marc C. and O'Hagan, Anthony},

  title = {{Bayesian} calibration of computer models},
  journal = {Journal of the Royal Statistical Society: Series B (Statistical Methodology)},
  year    = {2001},
  volume  = {63},
  number  = {3},
  pages   = {425--464},
}

@article{xu2024renewableClimate,

  author = {Xu, Luo and Feng, Kairui and Lin, Ning and Perera, A. T. D. and Poor, H. Vincent and Xie, Le and Ji, Chuanyi and Sun, X. Andy and Guo, Qinglai and O'Malley, Mark},

  title = {Resilience of renewable power systems under climate risks},
  journal = {Nature Reviews Electrical Engineering},
  year    = {2024},
  volume  = {1},
  pages   = {53--66},
}

@article{ruan2024interdisciplinary,

  author = {Ruan, Jiaqi and Xu, Zhao and Su, Hui},

  title = {Towards interdisciplinary integration of electrical engineering and earth science},
  journal = {Nature Reviews Electrical Engineering},
  year    = {2024},
  volume  = {1},
  pages   = {278--279},
}

@article{ruan2026longHorizon,

  author = {Ruan, Jiaqi and Li, Yishao and Gu, Chenjia and Jia, Youwei and Wang, Jianxiao and Hu, Zhengyang and Xu, Zhao and Su, Hui and Lie, Tek Tjing and Song, Jie},

  title = {Long-horizon resilience of power systems under climate change},
  journal = {Nature Reviews Electrical Engineering},
  year    = {2026},
  volume = {3},
  pages = {651--662},
}

@article{zscheischler2018compound,

  author = {Zscheischler, Jakob and Westra, Seth and van den Hurk, Bart J. J. M. and Seneviratne, Sonia I. and Ward, Philip J. and Pitman, Andy and AghaKouchak, Amir and Bresch, David N. and Leonard, Michael and Wahl, Thomas and Zhang, Xuebin},

  title = {Future climate risk from compound events},
  journal = {Nature Climate Change},
  year    = {2018},
  volume  = {8},
  pages   = {469--477},
}

@article{zscheischler2020typology,

  author = {Zscheischler, Jakob and Martius, Olivia and Westra, Seth and Bevacqua, Emanuele and Raymond, Colin and Horton, Radley M. and van den Hurk, Bart and AghaKouchak, Amir and J{\'e}z{\'e}quel, Agla{\'e} and Mahecha, Miguel D. and Maraun, Douglas and Ramos, Alexandre M. and Ridder, Nina N. and Thiery, Wim and Vignotto, Edoardo},

  title = {A typology of compound weather and climate events},
  journal = {Nature Reviews Earth \& Environment},
  year    = {2020},
  volume  = {1},
  pages   = {333--347},
}

@article{raymond2020connected,

  author = {Raymond, Colin and Horton, Radley M. and Zscheischler, Jakob and Martius, Olivia and AghaKouchak, Amir and Balch, Jennifer and Bowen, Steven G. and Camargo, Suzana J. and Hess, Jeremy and Kornhuber, Kai and Oppenheimer, Michael and Ruane, Alex C. and Wahl, Thomas and White, Kathleen},

  title = {Understanding and managing connected extreme events},
  journal = {Nature Climate Change},
  year    = {2020},
  volume  = {10},
  pages   = {611--621},
}

@article{surminski2016climateResilience,

  author = {Surminski, Swenja and Bouwer, Laurens M. and Linnerooth-Bayer, Joanne},

  title = {How insurance can support climate resilience},
  journal = {Nature Climate Change},
  year    = {2016},
  volume  = {6},
  pages   = {333--334},
}

@article{venkatasubramanian2023resilienceBibliometric,

  author = {Venkatasubramanian, Balaji V. and Panteli, Mathaios},

  title = {Power system resilience during 2001--2022: a bibliometric and correlation analysis},
  journal = {Renewable and Sustainable Energy Reviews},
  year    = {2023},
  volume  = {188},
  pages   = {113862},
}

@article{zidane2025microgridResilience,

  author = {Zidane, Tekai Eddine Khalil and Muis, Zarina Ab and Ho, Wai Shin and Zahraoui, Younes and Aziz, Ali Saleh and Su, Chun-Lien and Mekhilef, Saad and Campana, Pietro Elia},

  title = {Power systems and microgrids resilience enhancement strategies: a review},
  journal = {Renewable and Sustainable Energy Reviews},
  year    = {2025},
  volume  = {207},
  pages   = {114953},
}

@article{he2024cyberReview,

  author = {He, Rong and Jin, Zhuo and Li, Johnny Siu-Hang},

  title = {Modeling and management of cyber risk: a cross-disciplinary review},
  journal = {Annals of Actuarial Science},
  year    = {2024},
  volume  = {18},
  number  = {2},
  pages   = {270--309},
}

@article{eling2022cyberSeverity,

  author = {Eling, Martin and Jung, Kwangmin},

  title = {Heterogeneity in cyber loss severity and its impact on cyber risk measurement},
  journal = {Risk Management},
  year    = {2022},
  volume  = {24},
 number  = {4},
  pages   = {273--297},
}

@article{eling2022cyberQuantile,

  author = {Eling, Martin and Jung, Kwangmin and Shim, Jeungbo},

  title = {Unraveling heterogeneity in cyber risks using quantile regressions},
  journal = {Insurance: Mathematics and Economics},
  year    = {2022},
  volume  = {104},
  pages   = {222--242},
}

@article{hudson2016insuranceAdaptation,

  author = {Hudson, Paul and Botzen, W. J. Wouter and Feyen, Luc and Aerts, Jeroen C. J. H.},

  title = {Incentivising flood risk adaptation through risk based insurance premiums: trade-offs between affordability and risk reduction},
  journal = {Ecological Economics},
  year    = {2016},
  volume  = {125},
  pages   = {1--13},
}

@article{mol2020moralHazard,

  author = {Mol, Jantsje M. and Botzen, W. J. Wouter and Blasch, Julia E.},

  title = {Risk reduction in compulsory disaster insurance: experimental evidence on moral hazard and financial incentives},
  journal = {Journal of Behavioral and Experimental Economics},
  year    = {2020},
  volume  = {84},
  pages   = {101500},
}

@article{kousky2021parametric,

  author = {Kousky, Carolyn and Wiley, Helen and Shabman, Len},

  title = {Can parametric microinsurance improve the financial resilience of low-income households in the {United States}? {A} proof-of-concept examination},
  journal = {Economics of Disasters and Climate Change},
  year    = {2021},
  volume  = {5},
  pages   = {301--327},
}

@article{steptoe2022catModeling,

  author = {Steptoe, Hamish and Souch, Claire and Slingo, Julia},

  title = {Advances in numerical weather prediction, data science, and open-source software herald a paradigm shift in catastrophe risk modeling and insurance underwriting},
  journal = {Risk Management and Insurance Review},
  year    = {2022},
  volume  = {25},
  number  = {1},
  pages   = {69--81},
}

@article{wu2020catInsuranceEquilibrium,

  author = {Wu, Yang-Che},

  title = {Equilibrium in natural catastrophe insurance market under disaster-resistant technologies, financial innovations and government interventions},
  journal = {Insurance: Mathematics and Economics},
  year    = {2020},
  volume  = {95},
  pages   = {116--128},
}

@techreport{cebotari2020sovereignInsurance,

  author = {Cebotari, Aliona and Youssef, Karim},
  title       = {Natural Disaster Insurance for Sovereigns: issues, Challenges and Optimality},
  institution = {International Monetary Fund},
  type        = {IMF Working Paper},
  number      = {2020/003},
  year        = {2020},
}

@article{bracken2025futureDrought,

  author = {Bracken, Cameron and Voisin, Nathalie and Mongird, Kendall and Burleyson, Casey D. and Oikonomou, Konstantinos},

  title = {Intensifying renewable energy droughts in the western {U.S.} amid evolving infrastructure and climate},
  journal = {Earth's Future},
  year    = {2025},
  volume  = {13},
  number  = {6},
  pages   = {e2024EF005313},
}

@article{xu2026insuranceBudget,

  title = {Resilience enhancement strategy for distribution system considering insurance investment budgets},
  author = {Xu, Mingqian and Li, Gengfeng and Bian, Yiheng and Li, Minghao and Chen, Jun and Zhao, Tianyang and Bie, Zhaohong},
  journal = {IET Generation, Transmission \& Distribution},
  volume  = {20},
  number  = {1},
  pages   = {e70368},
  year    = {2026},
}

@article{wang2026hierarchicalInsurance,

  title = {Data driven resilience oriented hierarchical dispatch and insurance strategy for power systems under extreme weather events},
  author = {Wang, Xiuli and Li, Yonghao and Ullah, Zia and Wen, Fushuan and Wang, Shunli},
  journal = {Sustainable Energy, Grids and Networks},
  volume  = {47},
  pages   = {102382},
  year    = {2026},
}

@book{grossi2005catastropheModeling,
  title     = {Catastrophe Modeling: A New Approach to Managing Risk},

  publisher = {Springer New York, NY},
  year      = {2005},

  editor = {Grossi, Patricia and Kunreuther, Howard},
}

@article{tsanakas2005riskPricing,

  title = {Measurement and pricing of risk in insurance markets},
  author = {Tsanakas, Andreas and Desli, Evangelia},
  journal = {Risk Analysis},
  volume  = {25},
  number  = {6},
  pages   = {1653--1668},
  year    = {2005},
}

@article{huang2026insurerV2G,

  title = {{Enhancing power distribution system resilience through insurance mechanisms: An insurer optimization approach with V2G aggregator participation}},
  author = {Huang, Jing and Liao, Kai and Gu, Chenjia and Yang, Jianwei and Xiang, Yue and Shao, Chengcheng and Zhao, Alexis Pengfei},
  journal = {IEEE Transactions on Smart Grid},
volume  = {17},
  number  = {6},
  pages   = {6147--6161},
  year    = {2026},
}

@article{bracken2025seasonalDrought,

  author = {Bracken, Cameron and Voisin, Nathalie and Son, Youngjun and Feng, Sha and Anderson, Osten and Chen, Xiaodong and Li, He and Oikonomou, Konstantinos},

  title = {Seasonal compound renewable energy droughts in the {United States}},
  journal = {Environmental Research: Energy},
  year    = {2025},
  volume  = {2},
number  = {2},
  pages   = {025005},
}

@article{bian2025powerDrought,

  author = {Bian, Xiaoyan and Wang, Xueer and Zhou, Bo and Zhang, Jiawei and Wang, Tingting and Lin, Shunfu},

  title = {Spatial-temporal analysis of `power drought' under compound dry-hot events for renewable power systems},
  journal = {IET Renewable Power Generation},
  year    = {2025},
  volume  = {19},
  number  = {1},
  pages   = {e70120},
}

@incollection{chen2021ipccrisk,
  author    = {Chen, D. and Rojas, M. and Samset, B. H. and Cobb, K. and
               Diongue Niang, A. and Edwards, P. and Emori, S. and
               Faria, S. H. and Hawkins, E. and Hope, P. and Huybrechts, P. and
               Meinshausen, M. and Mustafa, S. K. and Plattner, G.-K. and
               Tr{\'e}guier, A.-M.},
  title     = {Framing, Context, and Methods},
  booktitle = {Climate Change 2021: The Physical Science Basis.
               Contribution of Working Group {I} to the Sixth Assessment
               Report of the Intergovernmental Panel on Climate Change},
  publisher = {Cambridge University Press},
  year      = {2021},
  pages     = {147--286},
  doi       = {10.1017/9781009157896.003},
  url       = {https://www.ipcc.ch/report/ar6/wg1/chapter/chapter-1/}
}

\end{document}